\documentclass{article}
\usepackage{latexsym}
\usepackage{amsfonts}
\usepackage{amssymb}
\usepackage{theorem}

\DeclareSymbolFont{rsfscript}{OMS}{rsfs}{m}{n}
\DeclareSymbolFontAlphabet{\mathrsfs}{rsfscript}

\theorembodyfont{\itshape}
\newtheorem{theorem}{Theorem}[section]
\newtheorem{corollary}{Corollary}[section]
\newtheorem{lemma}{Lemma}[section]
\newtheorem{proposition}{Proposition}[section]
\theorembodyfont{\rmfamily}

\newtheorem{remark}{Remark}[section]

\makeatletter

\@addtoreset{equation}{section}
\makeatother

\newcommand{\qed}{\hfill$\Box$}

\begin{document}

\title{Global wellposedness of the modified Benjamin-Ono equation with initial data in $H^{1/2}$}
\author{Carlos E. Kenig and Hideo Takaoka}
\date{\empty}
\maketitle
\begin{abstract}
We prove that the modified Benjamin-Ono equation is globally wellposed in $H^s$ for $s\ge 1/2$.
The exponent $H^{1/2}$ seems to be optimal in the sense that the solution map is not $C^3$ in $H^s$ for $s<1/2$ \cite{MR1}.
We perform a gauge transformation as in T. Tao \cite{Ta2}, but we combine it with a Littlewood-Paley decomposition.
We also use a space-time $L^2$-estimate that it is able to handle solutions in $H^{1/2}$ instead of solutions in the Besov space $B_{2,1}^{1/2}$ \cite{MR1}.
\end{abstract}
\section{Introduction}
In this paper, we consider the initial value problem for the modified Benjamin-Ono equation of the form\footnote{Also the equation with the nonlinearity of the form $-u^2u_x$ can be treated by our method.}
\begin{eqnarray}\label{eq:gBO}
\left\{
\begin{array}{l}
u_t+{\cal H}u_{xx}+u^2u_x=0,\quad (x,t)\in\Bbb R^2,\\
u(x,0)=u_0(x),
\end{array}
\right.
\end{eqnarray}
where $u:\Bbb R^2\to\Bbb R$ is a real-valued function and ${\cal H}$ is the Hilbert transform
\begin{eqnarray*}
{\cal H}u(x)=\frac{1}{\pi}\mbox{p.v.}\int_{-\infty}^{+\infty}\frac{u(y)}{x-y}\,dy.
\end{eqnarray*}
For the equation with quadratic nonlinearity
\begin{eqnarray}\label{eq:BO}
u_t+{\cal H}u_{xx}+(u^2)_x=0
\end{eqnarray}
Benjamin \cite{Be} and Ono \cite{O} derived this as a model for one-dimensional waves in deep water.
On the other hand, the cubic nonlinearity, found in a manner analogous to the relation between the KdV equation and the modified KdV equation, is also of much interest for long wave models, \cite{ABFS,KPV1}.

Recall that the conservation laws provide a priori bounds on the solution; namely there are at least the following three conservation laws preserved under the flow\footnote{For the Benjamin-Ono equation, the equation is completely integrable, and in fact possesses an infinite number of conservation laws.}
\begin{eqnarray}
& & \frac{d}{dt}\int_{\Bbb R}u(x,t)\,dx=0\nonumber\\
& & \frac{d}{dt}\int_{\Bbb R}u^2(x,t)\,dx=0\quad\mbox{($L^2$-mass)},\label{eq:con1}\\
& & \frac{d}{dt}\int_{\Bbb R}\frac{1}{2}u{\cal H}u_x-\frac{1}{12}u^{4}\,dx=0\quad\mbox{(Hamiltonian)}.\label{eq:con2}
\end{eqnarray}
Then establishing a global solution on the Hilbert space $H^{1/2}$ is of interest by the preservation of the Hamiltonian and the $L^2$-mass.
The purpose of this paper is in particular to prove the global wellposedness for data $u_0\in H^s$, for $s\ge 1/2$.
\begin{theorem}\label{thm:k=2}
Let $s\ge 1/2$.
For any $u_0\in H^{s}$, there exist $T=T(\|u_0\|_{H^{1/2}})$ and a unique solution $u$ of the equation (\ref{eq:gBO}) satisfying
\begin{eqnarray*}
u\in C([-T,T]:H^s)\cap X_{T}^s,
\end{eqnarray*}
where we shall define later the function space $X_{T}^s$ (see the end of this section).
Moreover, for any $R>0$ the solution operator $u_0\mapsto u$ is Lipschitz continuous from $\{u_0\in H^{s}:\|u_0\|_{H^s}<R\}$ to $C([-T,T]:H^{s})$.
\end{theorem}
We make some remarks about Theorem \ref{thm:k=2}.
\begin{remark}\label{rem:ill-posedness}
(i) Recall that heuristically the scaling argument
\begin{eqnarray}\label{eq:scaling}
u(x,t)\mapsto u_{\lambda}=\frac{1}{\lambda^{1/2}}u(\frac{x}{\lambda},\frac{t}{\lambda^2})
\end{eqnarray}
leads the constraint $s\ge 0$ on the wellposedness for (\ref{eq:gBO}).
The result in Theorem \ref{thm:k=2} is far from those given by scaling.
\newline
(ii) It is worth noting that when $s<1/2$, the solution map $u_0\mapsto u$ as mapping from $H^s$ to $C([-T,T]:H^s)$ is no longer of class $C^3$ \cite{MR1}.
Note that this illposedness result is true not only for $H^s$ but also for $B_{2,1}^s$.
Thus the value of $s=1/2$ in Theorem \ref{thm:k=2} may relate to the lower threshold of the result on local wellposedness. (It is not clear whether the solution map, given by Theorem \ref{thm:k=2}, is of $C^3$-class or not.)
\end{remark}

From the conservation laws (\ref{eq:con1})-(\ref{eq:con2}), and iterating Theorem \ref{thm:k=2} we obtain the following corollary.
\begin{corollary}
The Cauchy problem (\ref{eq:gBO}) is globally wellposed in $H^s$ for $s\ge 1/2$.
\end{corollary}

The initial value problem for the equation (\ref{eq:gBO}) and for the Benjamin-Ono equation (\ref{eq:BO}) have been extensively studied \cite{BL,CKS,GV1,GV2,GV3,I,KPV3,KK,KT,MR1,MR2,P,Ta2}; for instance the energy method provides the wellposedness on the Sobolev space $H^s$ for $s>3/2$ (see \cite{I}).
For the Benjamin-Ono equation (\ref{eq:BO}), it has been known that this is locally wellposed for $s>9/8$ \cite{KK,KT} by the refinement of the energy method and dispersive estimates.
T. Tao \cite{Ta2} extended this result to the energy space $s\ge 1$.
More precisely, the global wellposedness was obtained from the conservation law
\begin{eqnarray*}
\frac{d}{dt}\int_{\Bbb R}\frac12u_x^2+\frac32u^2{\cal H}u_x-\frac14u^4\,dx=0
\end{eqnarray*}
and the use of a "gauge transformation", where the solution map $u_0\mapsto u(t)$ persists in the $H^s$, but the Lipschitz continuity holds only in the $L^2$-space.
Very recently, the $H^1$-result was improved by A. D. Ionescu and C. E. Kenig \cite{IK} which obtained global wellposedness for $s\ge 0$, and also by N. Burq and F. Planchon \cite{BP} which obtained local wellposedness for $s>1/4$.  

For the modified Benjamin equation (\ref{eq:gBO}), L. Molinet and F. Ribaud \cite{MR2} have shown the local wellposedness in the Sobolev space $H^s$ for $s>1/2$.
(Results for the generalized Benjamin-Ono equation with higher nonlinearities, are also found in \cite{MR2}).
Their proof is based on Tao's gauge transformation.
Also, the result for $s=1/2$, but with the Sobolev space $H^s$ replaced by the Besov space $B_{2,1}^s$, has been obtained in \cite{MR1}; more precisely, they have proved local wellposedness in $B_{2,1}^s$ for $s\ge 1/2$.
In this result, however the smallness condition on the data is required.

Our method relies on a refinement of the gauge transformation (using a Littlewood-Paley decomposition), introduced initially for the Benjamin-Ono equation \cite{Ta2} and modified for the generalized Benjamin-Ono equation \cite{MR2}, as well as the use of estimates for the Schr\"odinger equation.
The point is that we shall transform the equation (\ref{eq:gBO}) into a derivative nonlinear Schr\"odinger equation, where the nonlinearity $u^2u_x$ in (\ref{eq:gBO}) has been placed relatively close to the form $\sum_{N_{high}}\partial_x(\sum_{N_{low}\ll N_{high}}P_{N_{low}}u)^2P_{N_{high}}u$, in other words, the derivative in the nonlinearity does not appear in the highest frequency terms.
We will describe this reduction of the equation in the next section.
\begin{remark}\label{lem:DNLS}
A very similar equation to (\ref{eq:gBO}) is the derivative nonlinear Schr\"odinger equation
\begin{eqnarray}\label{eq:DNLS}
u_t-iu_{xx}+|u|^2u_x=0,~(x,t)\in\Bbb R^2,
\end{eqnarray}
and local wellposedness was known for the equation in $H^s$ for $s\ge 1/2$ \cite{Tak}, where a fixed point argument is performed in an adapted Bourgain's $X_{s,b}$ space which yields a $C^{\infty}$-solution map.
Our method also gives the result for the equation (\ref{eq:DNLS}), without Bourgain's space, in $H^s$ for $s\ge 1/2$, but only shows the solution map to be Lipshitz.
\end{remark}

One difficulty in proving the ``endpoint'' case $s=1/2$ for solutions in $H^s$, is that using the inhomogeneous smoothing effect estimate (see c.f. (\ref{eq:rsmoothing-2}) below)
\begin{eqnarray*}
\left\|\partial_x\int_0^te^{-t(t-t'){\cal H}\partial_x^2}f(t')\,dt'\right\|_{L_x^{\infty}L_T^2}\lesssim\|f\|_{L_x^1L_T^2},
\end{eqnarray*}
one loses the $L_x^1L_T^2$ control for the nonlinearity $u^2u_x$.
This is because one needs to use the $L_x^2L_T^{\infty}$-maximal function estimate for two $u$' and the $L_x^{\infty}L_T^2$-smoothing effect estimate for the term $u_x$, and in principle the maximal function estimate may fail at the endpoint $s=1/2$, although the estimate is valid at the endpoint provided that the data are dyadically localized in frequency space.
In fact, we use the $l^2$-type maximal function estimate in order to invoke the endpoint maximal function estimate
\begin{eqnarray*}
\left(\sum_N\|e^{-t{\cal H}\partial_x^2}P_Nu_0\|_{L_x^2L_T^{\infty}}^2\right)^{1/2}\lesssim \|u_0\|_{H^{1/2}}.
\end{eqnarray*}
We then estimate the $L_x^1L_T^2$-type norm for the nonlinearity.
To summarize, we suppose the nonlinearity to be $\sum_{N_{high}}\partial_x(\sum_{N_{low}\ll N_{high}}P_{N_{low}}u)^2P_{N_{high}}u$ as mentioned before.
Applying the Littlewood-Paley projection operator $P_N$ to the equation, for each $N$, we estimate this by
\begin{eqnarray*}
& & \left\|P_N\left(\sum_{N_{high}}\partial_x(\sum_{N_{low}\ll N_{high}}P_{N_{low}}u)^2P_{N_{high}}u\right)\right\|_{L_x^1L_T^2}\\
 & \lesssim & \sum_{N_{high}\sim N}\left\|\partial_x(\sum_{N_{low}\ll N_{high}}P_{N_{low}}u)^2\right\|_{L_{xT}^2}\|P_{N_{high}}u\|_{L_x^2L_T^{\infty}}.
\end{eqnarray*}
In particular, we prove the following space-time $L^2$ estimate which is crucial to our proof of Theorem \ref{thm:k=2} (see section \ref{sec:L^2-spacetime} for the proof of this proposition).
\begin{proposition}\label{prop:L^2-spacetime}
Let $u$ be a $H^{\infty}$-solution to (\ref{eq:gBO}).
Then we have
\begin{eqnarray}\label{eq:L^2-spacetime}
\|(u^2)_x\|_{L_{xT}^2}\lesssim \|P_{\ge 1}u_0\|_{H^{1/2}}^2+T^{\frac12}\|u\|_{X_T^{1/2}}^2+(1+\|u\|_{X_T^{1/2}})\|u\|_{X_T^{1/2}}\|P_{\ge 1}u\|_{X_T^{1/2}}.
\end{eqnarray}
\end{proposition}
We close this section by introducing some notation.
Let $\psi$ be a fixed even $C^{\infty}$ function of compact support, with $\mbox{supp}\psi\subset\{|\xi|<2\}$, and $\psi(\xi)=1$ for $|\xi|\le 1$.
Define $\varphi(\xi)=\psi(\xi)-\psi(2\xi)$.
Let $N$ be a dyadic number of the form $N=2^j,~j\in\Bbb N\cup\{0\}$ or $N=0$.
Writing $\varphi_N(\xi)=\varphi(\xi/N)$ for $N\ge 1$, we define the convolution operator $P_N$ by $P_Nu=u*\check{\varphi}_N$, where $\check{\cdot}$ denotes spatial Fourier inverse transform (while $\hat{\cdot}$ denotes a spatial Fourier transform).
Then we have a spatial Littlewood-Paley decomposition
\begin{eqnarray*}
\sum_{N}P_N=I
\end{eqnarray*}
where we define the function $\varphi_0$ by $\varphi_0(\xi)=1-\sum_{N}\varphi_N(\xi)$ to denote $P_0u=u*\check{\varphi}_0$.
Note that if $u$ is real-valued function, then $P_Nu$ is also real-valued.
We define the projection operators $P_{\pm}$ to the frequency $\pm[0,\infty)$.
We will recall the Littlewood-Paley theorem \cite{S}
\begin{eqnarray*}
\left\|\left(\sum_N|P_N\phi|^2\right)^{1/2}\right\|_{L^p}\sim \|\phi\|_{L^p}
\end{eqnarray*}
for $1<p<\infty$.

For nonnegative quantities $A,B$, we use $A\lesssim B$ to denote the estimate $A\le CB$ for some $C>0$, and $A\sim B$ to denote $A\lesssim B\lesssim A$, and $A\ll B$ to denote $A\le CB$ for some small $C>0$.

We also define more general projection $P_{\ll N}$ and $P_{\lesssim N}$ by
\begin{eqnarray*}
P_{\ll N}=\sum_{M\ll N}P_M,~P_{\lesssim N}=\sum_{M\lesssim N}P_M.
\end{eqnarray*}
Similarly define $P_{\gg N},~P_{\gtrsim N}$ and also define $P_{\sim N}$, etc.
We remark that the projection operators $P_{\ll N},P_{\lesssim N},P_{\gg N},P_{\gtrsim N}$ are bounded on $L^p,L_T^qL_x^p,L_x^pL_T^q$, for $1<p,q<\infty$.
Moreover, $P_N$ and ${\cal H}P_N$ are bounded operators on $L^p,L_T^qL_x^p,L_x^pL_T^q$ for $1\le p,q\le\infty$.

We define the Lebesgue spaces $L_T^qL_x^p$ and $L_x^pL_T^q$ by the norms
\begin{eqnarray*}
\|f\|_{L_T^qL_x^p}=\|\|f\|_{L_x^p(\Bbb R)}\|_{L_t^q([0,T])},\quad
\|f\|_{L_x^pL_T^q}=\|\|f\|_{L_t^q([0,T])}\|_{L_x^p(\Bbb R)}.
\end{eqnarray*}
In particular, when $p=q$, we abbreviate $L_T^qL_x^p$ or $L_x^pL_T^q$ as $L_{xT}^p$.

Let $\langle\cdot\rangle=(1+|\cdot|^2)^{1/2}$.
We use the fractional differential operators $D_x^s$ and $\langle D_x\rangle^s$ defined by
\begin{eqnarray*}
\widehat{D_x^s f}(\xi)=|\xi|^s\widehat{f}(\xi),~\widehat{\langle D_x\rangle^s f}(\xi)=\langle \xi\rangle^s\widehat{f}(\xi).
\end{eqnarray*}

We are now ready to define the function space $X_{T}^s$.
For $s\ge 1/2,~T>0$, we introduce
\begin{eqnarray*}
X_{T}^s=\{u\in\mathrsfs{D}'(\Bbb R\times (-T,T)):\|u\|_{X_T^s}<\infty\},
\end{eqnarray*}
where
\begin{eqnarray*}
\|u\|_{X_{T}^s}& = & \|u\|_{L_T^{\infty}H^s}+\left(\sum_{k=1}^{[s+1/2]}\sum_{N}\|D_x^{s+\frac12-k}\partial_x^{k}P_Nu\|_{L_x^{\infty}L_T^2}^2\right)^{1/2}\\
& & +\left(\sum_{k=0}^{[s]}\sum_{N}\|\langle D_x\rangle^{s-k-\frac12}\partial_x^{k}P_Nu\|_{L_x^2L_T^{\infty}}^2\right)^{1/2}+\left(\sum_{k=0}^{[s]}\sum_{N}\|\langle D_x\rangle^{s-k-\frac14}\partial_x^{k}P_Nu\|_{L_x^4L_T^{\infty}}^2\right)^{1/2}.
\end{eqnarray*}

Acknowledgments.
This research was performed while the second author (H.T.) visited Carlos E. Kenig at the University of Chicago under the J.S.P.S. Postdoctoral Fellowships for Research Abroad.
He would like to thank the University of Chicago for its hospitality.
\section{Gauge transformation}
We transform the equation (\ref{eq:gBO}) as stated in the introduction.
Let $u(x,t)$ be an $H^{\infty}$-solution to (\ref{eq:gBO}).
We introduce the complex-valued functions $v_N:\Bbb R^2\to\Bbb C$ for a dyadic number $N$ by\footnote{This gauge transform is also inspired by the result in \cite{Tak}. In fact, when $u$ is a complex-valued function and the nonlinearity in (\ref{eq:gBO}) is replaced by $|u|^2u_x$, we let $v_N(x,t)=e^{-\frac{i}{2}\int^x|P_{\ll N}u|^{2}}P_+P_Nu(x,t)$; which is modified from $v(t,x)=e^{-\frac{i}{2}\int^x|u|^2}u(x,t)$ in \cite{Tak}. We also mention that if $v_N(x,t)=e^{-\frac{i}{2}\int^xP_{\ll N}u}P_+P_Nu(x,t)$, our method gives the $H^1$-wellposedness for the Benjamin-Ono equation \cite{Ta2}.}
\begin{eqnarray}\label{eq:gauge}
v_N(x,t)=e^{-\frac{i}{2}\int_{-\infty}^x(P_{\ll N}u(y,t))^2\,dy}P_+P_Nu(x,t).
\end{eqnarray}
It will be convenient to abbreviate by writing $e^{-\frac{i}{2}\int^x(P_{\ll N}u)^2}$ for $e^{-\frac{i}{2}\int_{-\infty}^x(P_{\ll N}u(y,t))^2\,dy}$.
From the Leibniz rule and ${\cal H}P_+=-i$, we see that
\begin{eqnarray*}
(\partial_t-i\partial_x^2)v_N &=& e^{-\frac{i}{2}\int^x(P_{\ll N}u)^2}(\partial_t+{\cal H}\partial_x^2)P_+P_Nu\\
& & -e^{-\frac{i}{2}\int^x(P_{\ll N}u)^2}\partial_xP_+P_Nu\partial_x\int_{-\infty}^x(P_{\ll N}u(y,t))^2\,dy\\
& & +P_+P_Nu(\partial_t-i\partial_x^2)e^{-\frac{i}{2}\int^x(P_{\ll N}u)^2}\\
&=& e^{-\frac{i}{2}\int^x(P_{\ll N}u)^2}(P_+P_N(u^2u_x)-(P_{\ll N}u)^2P_+P_Nu_x)\\
& & -\frac{i}{2}e^{-\frac{i}{2}\int^x(P_{\ll N}u)^2}P_+P_Nu(\partial_t-i\partial_x^2)\int_{-\infty}^x(P_{\ll N}u(y,t))^2\,dy\\
& & +\frac{i}{4}e^{-\frac{i}{2}\int^x(P_{\ll N}u)^2}P_+P_Nu(P_{\ll N}u)^{4}.
\end{eqnarray*}
In carrying out the computation for the second term, we use the equation (\ref{eq:gBO}) and integrate by parts.
Thus
\begin{eqnarray*}
& &(\partial_t-i\partial_x^2)\int_{-\infty}^x(P_{\ll N}u(y,t))^2\,dy\\
&=& 2\int_{-\infty}^xP_{\ll N}u\partial_tP_{\ll N}u\,dy-2iP_{\ll N}uP_{\ll N}u_x\\
&=& -2\int_{-\infty}^xP_{\ll N}uP_{\ll N}({\cal H}u_{xx}+u^2u_x)\,dy-2iP_{\ll N}uP_{\ll N}u_x\\
&=& 2i(i{\cal H}-1)P_{\ll N}u_xP_{\ll N}u+2\int_{-\infty}^x{\cal H}P_{\ll N}u_xP_{\ll N}u_x\,dy\\
& & -2\int_{-\infty}^xP_{\ll N}uP_{\ll N}(u^2u_x)\,dy.
\end{eqnarray*}
Hence, $v_N$ finally obeys the following differential equation
\begin{eqnarray}
(\partial_t-i\partial_x^2)v_N & = & e^{-\frac{i}{2}\int^x(P_{\ll N}u)^2}(P_+P_N(u^2u_x)-(P_{\ll N}u)^2P_+P_Nu_x)\nonumber\\
& & +e^{-\frac{i}{2}\int^x(P_{\ll N}u)^2}P_{\ll N}(i{\cal H}-1)u_xP_{\ll N}uP_+P_Nu\nonumber\\
& & -ie^{-\frac{i}{2}\int^x(P_{\ll N}u)^2}P_+P_Nu\int_{-\infty}^x{\cal H}P_{\ll N}u_xP_{\ll N}u_x\,dy\nonumber\\
& & +ie^{-\frac{i}{2}\int^x(P_{\ll N}u)^2}P_+P_Nu\int_{-\infty}^xP_{\ll N}(u^2u_x)P_{\ll N}u\,dy\nonumber\\
& & +\frac{i}{4}e^{-\frac{i}{2}\int^x(P_{\ll N}u)^2}(P_{\ll N}u)^4P_+P_Nu\nonumber\\
& \equiv &A_{1,N}(t)+A_{2,N}(t)+A_{3,N}(t)+A_{4,N}(t)+A_{5,N}(t).
\label{eq:mmBO5}
\end{eqnarray}
The desired a priori estimate for $u$ in (\ref{eq:gBO}) can be proven from the solutions $v_N$ in (\ref{eq:mmBO5}).
We prove this in section \ref{sec:apriori}.
\begin{remark}
As opposed to (\ref{eq:gBO}), for (\ref{eq:mmBO5}), the very worst type of nonlinearity as $(P_{N_{low}}u)^2P_{N_{high}}u_x$ with $|N_{low}|\ll |N_{high}|$, in which the derivative on one of three $u$'s can not be shared with the other two $u$'s, is almost absent.
This is a consequence of the formula; for instance we expand
\begin{eqnarray}
& & P_+P_N(u^2u_x)-(P_{\ll N}u)^2P_+P_Nu_x\nonumber\\
&=& P_N((P_{\ll N}u)^2P_+\widetilde{P}_Nu_x)-(P_{\ll N}u)^2P_+P_N\widetilde{P}_Nu_x\nonumber\\
& & +P_+P_N((u^2-(P_{\ll N}u)^2)u_x)\label{eq:absent}
\end{eqnarray}
for $N\gg 1$, where $\widetilde{P}_N=P_{N/2}+P_N+P_{2N}$.
One can think in particular of the first term in (\ref{eq:absent}) as $c(P_{\ll N}u)^2_xP_+\widetilde{P}_Nu$ (see section \ref{sec:nonlinear}).
\end{remark}
\section{Preliminaries}
In order to prove the a priori estimate for the equation of $v_N$, we need the linear estimates associated with the one-dimensional Schr\"odinger equation.
We first recall the Strichartz estimates, smoothing effects and maximal function estimates (for the proof, see e.g. \cite{KPV2}).
\begin{lemma}\label{lem:freeestimate}
For all $\phi\in\mathrsfs{S}(\Bbb R),~\theta\in[0,1]$ and $T\in(0,1)$,
\begin{eqnarray}\label{eq:strichartz}
\|e^{it\partial_x^2}\phi\|_{L_T^{\frac{4}{\theta}}L_x^{\frac{2}{1-\theta}}}\lesssim\|\phi\|_{L^2},
\end{eqnarray}
\begin{eqnarray}\label{eq:smoothing}
\|e^{it\partial_x^2}P_N\phi\|_{L_x^{\frac{2}{1-\theta}}L_T^{\frac{2}{\theta}}}\lesssim \langle N\rangle^{\frac{1}{2}-\theta}\|\phi\|_{L^2},
\end{eqnarray}
\begin{eqnarray}\label{eq:maximal}
\|e^{it\partial_x^2}\phi\|_{L_x^4L_T^{\infty}}\lesssim \|\phi\|_{\dot{H}^{1/4}}.
\end{eqnarray}
\end{lemma}
\begin{remark}
We say that a pair $(q,p)$ is admissible if $\frac{2}{q}=\frac12-\frac{1}{p}$.
Then the above pair $(\frac{4}{\theta},\frac{2}{1-\theta})$ is admissible.
\end{remark} 
{\it Proof.}
The inequalities (\ref{eq:strichartz}) and (\ref{eq:maximal}) are due to the standard Strichartz and maximal function estimates, respectively \cite{KPV2}.

To show (\ref{eq:smoothing}) we need the following inequalities \cite{KPV2}
\begin{eqnarray}
& & \|\langle D_x\rangle^{\frac{1}{2}+i\alpha}e^{it\partial_x^2}\phi\|_{L_x^{\infty}L_T^2}\lesssim\|\phi\|_{L^2}\label{eq:smoothing-1},\\
& & \|\langle D_x\rangle^{i\alpha}e^{it\partial_x^2}P_N\phi\|_{L_x^2L_T^{\infty}}\lesssim \langle N\rangle^{1/2}\|\phi\|_{L^2}.\nonumber
\end{eqnarray}
Applying complex interpolation argument to these inequalities, we obtain (\ref{eq:smoothing}) (if necessary, we use the trivial inequality $\|e^{it\partial_x^2}P_0\phi\|_{L_{xT}^{\infty}}\lesssim\|\phi\|_{L^2}$ to justify (\ref{eq:smoothing-1})).
\qed

We next state the $L_T^qL_x^p$ and $L_x^pL_T^q$ estimates for the linear operator $f\mapsto\int_0^te^{i(t-t')\partial_x^2}f(t')\,dt'$.
\begin{lemma}\label{lem:retard}
For $f\in\mathrsfs{S}(\Bbb R^2),\theta\in[0,1]$ and $T\in(0,1)$,
\begin{eqnarray}\label{eq:rstrichartz-1}
\left\|\int_0^te^{i(t-t')\partial_x^2}f(t')\,dt\right\|_{L_T^{4/\theta}L_x^{2/(1-\theta)}}\lesssim \|f\|_{L_T^{(4/\theta)'}L_x^{(2/(1-\theta))'}}.
\end{eqnarray}
\begin{eqnarray}\label{eq:rstrichartz-2}
\left\|\langle D_x\rangle^{\frac{\theta}{2}}\int_0^te^{i(t-t')\partial_x^2}f(t')\,dt'\right\|_{L_T^{\infty}L_x^2}\lesssim \|f\|_{L_x^{p(\theta)}L_T^{q(\theta)}},
\end{eqnarray}
\begin{eqnarray}\label{eq:rsmoothing}
\left\|D_x^{\frac{1+\theta}{2}}\int_0^te^{i(t-t')\partial_x^2}f(t')\,dt'\right\|_{L_x^{\infty}L_T^{2}}\lesssim\|f\|_{L_x^{p(\theta)}L_T^{q(\theta)}},
\end{eqnarray}
\begin{eqnarray}\label{eq:rmaximal-1}
\left\|\langle D_x\rangle^{\frac{\theta}{2}}\int_0^te^{i(t-t')\partial_x^2}P_Nf(t')\,dt'\right\|_{L_x^{2}L_T^{\infty}}\lesssim \langle N\rangle^{1/2}\|f\|_{L_x^{p(\theta)}L_T^{q(\theta)}},
\end{eqnarray}
\begin{eqnarray}\label{eq:rmaximal-2}
\left\|\langle D_x\rangle^{\frac{\theta}{2}-\frac14}\int_0^te^{i(t-t')\partial_x^2}f(t')\,dt'\right\|_{L_T^{4}L_x^{\infty}}\lesssim\|f\|_{L_x^{p(\theta)}L_T^{q(\theta)}},
\end{eqnarray}
where $p'$ of number is conjugate of $p\in[1,\infty]$ given by $1/p+1/p'=1$, and
\begin{eqnarray*}
\frac{1}{p(\theta)}=\frac{3+\theta}{4},~
\frac{1}{q(\theta)}=\frac{3-\theta}{4}.
\end{eqnarray*}
\end{lemma}

We shall need the lemma of Christ-Kiselev \cite{CK}, which permits us to obtain Lemma \ref{lem:retard} from the corresponding "non-retarded estimates" (see also \cite{MR1,MR2,SS,Ta3}).
\begin{lemma}[Christ-Kiselev \cite{CK}]\label{lem:ChristKiselev}
Let $T$ be a linear operator of the form
\begin{eqnarray*}
Tf(t)=\int_{-\infty}^{\infty}K(t,t')f(t')\,dt'
\end{eqnarray*}
where $K:\mathrsfs{S}(\Bbb R^2)\to C(\Bbb R^3)$.
Assume that $\|Tf\|_{L_x^{p_1}L_T^{q_1}}\lesssim\|f\|_{L_x^{p_2}L_T^{q_2}}$ for $p_1,p_2,q_1,q_2\in[1,\infty]$ with $\min\{p_1,q_1\}>\max\{p_2,q_2\}$ or $p_2,q_2<\infty,q_1=\infty$.
Then
\begin{eqnarray*}
\left\|\int_0^tK(t,t')f(t')\,dt'\right\|_{L_x^{p_1}L_T^{q_1}}\lesssim\|f\|_{L_x^{p_2}L_T^{q_2}}.
\end{eqnarray*}
\end{lemma}
\begin{remark}\label{rem:ChristKiselev}
The $L_T^{q_1}L_x^{p_1},L_x^{p_2}L_T^{q_2}$ (instead of $L_T^{q_2}L_x^{p_2}$) version of Lemma \ref{lem:ChristKiselev} holds with the condition $q_1>\max\{p_2,q_2\}$ \cite{MR1,MR2}.
\end{remark}
{\it Proof of Lemma \ref{lem:retard}.}
The inequality (\ref{eq:rstrichartz-1}) is due to the inhomogeneous Strichartz estimate \cite{KPV2}.

The inequality (\ref{eq:rstrichartz-2}) follows from a $TT^*$ argument, (\ref{eq:smoothing-1}) and (\ref{eq:rstrichartz-1}).
Indeed, applying a $TT^*$ argument to (\ref{eq:smoothing-1}) we have
\begin{eqnarray*}
\left\|\langle D_x\rangle^{1/2}\int_{-\infty}^{\infty}e^{i(t-t')\partial_x^2}f(t')\,dt'\right\|_{L_T^{\infty}L_x^2}\lesssim \|f\|_{L_x^1L_T^2}.
\end{eqnarray*}
Also by (\ref{eq:rstrichartz-1}) we have
\begin{eqnarray*}
\left\|\int_{0}^{t}e^{i(t-t')\partial_x^2}f(t')\,dt'\right\|_{L_T^{\infty}L_x^2} & \lesssim & \min\{\|f\|_{L_T^1L_x^2},\|f\|_{L_T^{4/3}L_x^1}\}\\
& \lesssim & \min\{\|f\|_{L_{T}^{4/3}L_x^2},\|f\|_{L_{T}^{4/3}L_x^1}\}.
\end{eqnarray*}
Therefore by Remark \ref{rem:ChristKiselev} and the complex interpolation argument, we obtain (\ref{eq:rstrichartz-2}).

For (\ref{eq:rsmoothing}), in analogy with (\ref{eq:rstrichartz-2}) we begin with the following estimate
\begin{eqnarray}\label{eq:TT}
\left\|D_x^{\frac{\theta}{2}}\int_{-\infty}^{\infty}e^{i(t-t')\partial_x^2}f(t')\,dt'\right\|_{L_x^2}\lesssim\|f\|_{L_x^{p(\theta)}L_T^{q(\theta)}}.
\end{eqnarray}
This follows easily from the argument as before.
Then we use again $TT^*$ argument, (\ref{eq:smoothing-1}) and (\ref{eq:TT}) to obtain
\begin{eqnarray*}
\left\|D_x^{\frac{1+\theta}{2}}\int_{-\infty}^{\infty}e^{i(t-t')\partial_x^2}f(t')\,dt'\right\|_{L_x^{\infty}L_T^{2}} & \lesssim & \left\|D_x^{\frac{\theta}{2}}\int_{-\infty}^{\infty}e^{-it'\partial_x^2}f(t')\,dt'\right\|_{L_x^2}\\
& \lesssim & \|f\|_{L_x^{p(\theta)}L_T^{q(\theta)}}.
\end{eqnarray*}
Thus Lemma \ref{lem:ChristKiselev} implies (\ref{eq:rsmoothing}).

The proofs for (\ref{eq:rmaximal-1}) and (\ref{eq:rmaximal-2}) are the same as that for (\ref{eq:rsmoothing}) by using (\ref{eq:smoothing}) and (\ref{eq:maximal}).
\qed
\begin{remark}\label{rem:sretard}
(i) A straightforward application of Lemma \ref{lem:ChristKiselev} to (\ref{eq:smoothing}), (\ref{eq:maximal}), (\ref{eq:smoothing-1}) shows that for $\theta\in[0,1]$
\begin{eqnarray}\label{eq:srsmoothing}
\left\|\langle D_x\rangle^{1/2}\int_0^te^{i(t-t')\partial_x^2}f(t')\,dt'\right\|_{L_x^{\infty}L_T^2}\lesssim\|f\|_{L_T^1L_x^2},
\end{eqnarray}
\begin{eqnarray}\label{eq:srmaximal-1}
\left\|\int_0^te^{i(t-t')\partial_x^2}P_Nf(t')\,dt'\right\|_{L_x^{\frac{2}{\theta}}L_T^{\frac{2}{1-\theta}}}\lesssim\langle N\rangle^{1/2-\theta}\|f\|_{L_T^1L_x^2},
\end{eqnarray}
\begin{eqnarray}\label{eq:srmaximal-2}
\left\|\int_0^te^{i(t-t')\partial_x^2}f(t')\,dt'\right\|_{L_x^{4}L_T^{\infty}}\lesssim\|f\|_{L_T^1\dot{H}_x^{1/4}}.
\end{eqnarray}
\newline
(ii) The estimate (\ref{eq:rsmoothing}) with $\theta=1$, but with the $D_x$-derivative replaced by $\partial_x$, still holds \cite{KPV2}
\begin{eqnarray}\label{eq:rsmoothing-2}
\left\|\partial_x\int_0^te^{i(t-t')\partial_x^2}f(t')\,dt'\right\|_{L_x^{\infty}L_T^2}\lesssim\|f\|_{L_x^1L_T^2}.
\end{eqnarray}
\end{remark}

The proof of the estimates with the regularity $s$ for $s\in[1/2,1)$ requires that we use the Leibniz' type rule with the fractional-order differentiation.
The first lemma will provide the Leibniz' rule for the bilinear form $fg$.
\begin{lemma}\label{lem:Leibnitz1}
Let $\alpha\in(0,1),\alpha_1,\alpha_2\in[0,\alpha],~p,p_1,p_2,q,q_1,q_2\in(1,\infty)$ with $\alpha=\alpha_1+\alpha_2,~\frac{1}{p}=\frac{1}{p_1}+\frac{1}{p_2},~\frac{1}{q}=\frac{1}{q_1}+\frac{1}{q_2}$.
Then
\begin{eqnarray*}
\|D_x^{\alpha}(fg)-D_x^{\alpha}fg-fD_x^{\alpha}g\|_{L_x^pL_T^q}\lesssim\|D_x^{\alpha_1}f\|_{L_x^{p_1}L_T^{q_1}}\|D_x^{\alpha_2}g\|_{L_x^{p_2}L_T^{q_2}}.
\end{eqnarray*}
Moreover, the case $q_1=\infty$ is allowed if $\alpha_1=0$.
Added to this, the case $(p,q)=(1,2)$ is also allowed.
\end{lemma}
{\it Proof.}
See \cite[Theorems A.8 and A.13]{KPV1}.

Next, we shall have the Leibniz' rule for a product of the form $e^{iF}g$ where $F$ is the spatial primitive of some function $f$.
\begin{lemma}\label{lem:Leibnitz2}
Let $\alpha\in(0,1),~p,p_1,p_2,q,q_1\in(1,\infty),q_2\in(0,\infty]$ with $\frac{1}{p}=\frac{1}{p_1}+\frac{1}{p_2},~\frac{1}{q}=\frac{1}{q_1}+\frac{1}{q_2}$, and let $F(x,t)=\int_{-\infty}^xf(y,t)\,dy$, with real-valued function $f$.
Then
\begin{eqnarray*}
\|D_x^{\alpha}(e^{iF}g)\|_{L_x^pL_T^q}\lesssim\|f\|_{L_x^{p_1}L_T^{q_1}}\|g\|_{L_x^{p_2}L_T^{q_2}}+\|\langle D_x\rangle^{\alpha}g\|_{L_x^pL_T^q}.
\end{eqnarray*}
\end{lemma}
{\it Proof.}
We write
\begin{eqnarray}\label{eq:Lebnitz2d}
e^{iF}=P_{0}e^{iF}+P_{\ge 1}e^{iF}.
\end{eqnarray}
For the first term in (\ref{eq:Lebnitz2d}), we easily obtain the bound by $\lesssim \|P_0e^{iF}\langle D_x\rangle^{\alpha}g\|_{L_x^pL_T^q}\le \|\langle D_x\rangle^{\alpha}g\|_{L_x^pL_T^q}$.
To estimate the second term in (\ref{eq:Lebnitz2d}), we apply Lemma \ref{lem:Leibnitz1} to obtain the bound
\begin{eqnarray*}
& \le & \|D_x^{\alpha}(P_{\ge 1}e^{iF})g\|_{L_x^pL_T^q}+\|P_{\ge 1}e^{iF}D_x^{\alpha}g\|_{L_x^pL_T^q}+c\|D_x^{\alpha}P_{\ge 1}e^{iF}\|_{L_x^{p_1}L_T^{q_1}}\|g\|_{L_x^{p_2}L_T^{q_2}}\\
& \lesssim & \|f\|_{L_x^{p_1}L_T^{q_1}}\|g\|_{L_x^{p_2}L_T^{q_2}}+\|D_x^{\alpha}g\|_{L_x^pL_T^q}.
\end{eqnarray*}
The estimate on the term $\|D_x^{\alpha}P_{\ge 1}e^{iF}\|_{L_x^{p_1}L_T^{q_1}}$ is clear, by adding an extra derivative $D_x^{1-\alpha}$, and the fact that the Hilbert transform operator ${\cal H}$ is bounded on $L_x^pL_T^q$ to itself, for $1<p,q<\infty$.
\qed

In order to control the integral type nonlinearity in (\ref{eq:mmBO5}), we need the following lemma.
\begin{lemma}\label{lem:Leibnitz3}
Let $\alpha,\alpha_1,\alpha_2\in[0,1],~p,p_1,p_2,q,q_1,q_2\in(1,\infty)$ with $\alpha+1=\alpha_1+\alpha_2,~\frac{1}{p}=\frac{1}{p_1}+\frac{1}{p_2},~\frac{1}{q}=\frac{1}{q_1}+\frac{1}{q_2}$.
Then
\begin{eqnarray*}
\|D_x^{\alpha}\int_{-\infty}^x{\cal H}f_xf_x\,dy\|_{L_x^pL_T^q}\lesssim\|D_x^{\alpha_1}f\|_{L_x^{p_1}L_T^{q_1}}\|D_x^{\alpha_2}f\|_{L_x^{p_2}L_T^{q_2}}.
\end{eqnarray*}
\end{lemma}
{\it Proof.}
See \cite[Lemma 6.1]{MR2}.
\section{Proof of Proposition \ref{prop:L^2-spacetime}}\label{sec:L^2-spacetime}
In this section we prove Proposition \ref{prop:L^2-spacetime}.
Throughout the section, we will use $\xi_{ij}$ to denote $\xi_i+\xi_j$, and also use $\xi_{ijk}$ etc.

Using a Littlewood-Paley decomposition, we write
\begin{eqnarray*}
u^2=\sum_{N_1,N_2}P_{N_1}uP_{N_2}u.
\end{eqnarray*}
We split the sum into three parts $N_1\sim N_2,~N_1\gg N_2,~N_1\ll N_2$.

In the treatment of the case $N_1\sim N_2$, we can share a derivative between $P_{N_1}u$ and $P_{N_2}u$.
In fact, by Plancherel's theorem and the inequality
\begin{eqnarray*}
\left(\sum_N\|\langle D_x\rangle^{1/2}P_Nu\|_{L_{xT}^4}^2\right)^{1/2}\lesssim
T^{1/4}\|u\|_{X_T^{1/2}}+\|P_{\ge 1}u\|_{X_T^{1/2}}
\end{eqnarray*}
(which follows by interpolation), we have the bound for this contribution to the left-hand side of (\ref{eq:L^2-spacetime}) by
\begin{eqnarray*}
\lesssim \sum_{N_1\sim N_2}\|\langle D_x\rangle^{\frac12}P_{N_1}u\|_{L_{xT}^4}\|\langle D_x\rangle^{\frac12}P_{N_2}u\|_{L_{xT}^4}
\lesssim T^{\frac12}\|u\|_{X_T^{1/2}}^2+\|P_{\ge 1}u\|_{X_T^{1/2}}^2.
\end{eqnarray*}

Next consider the case $N_1\gg N_2$ or $N_1\ll N_2$.
By symmetry, it will suffice to consider the case $N_1\ll N_2$.
If $N_1=O(1)$, the proof is easy.
In fact,
\begin{eqnarray*}
\sum_{N_1=O(1)\ll N_2}P_{N_1}uP_{N_2}u_x=\sum_{N_1=O(1)\ll N_2}\widetilde{P}_{N_2}(P_{N_1}uP_{N_2}u_x).
\end{eqnarray*}
With this and the Littlewood-Paley theorem, we have the bound for this contribution to the left-hand side of (\ref{eq:L^2-spacetime}) by
\begin{eqnarray*}
\lesssim\sum_{N_1=O(1)}\|P_{N_1}u\|_{L_x^2L_T^{\infty}}\left(\sum_{N_2\gg 1}\|P_{N_2}u_x\|_{L_x^{\infty}L_T^2}^2\right)^{\frac12}\lesssim \|u\|_{X_T^{1/2}}\|P_{\ge 1}u\|_{X_T^{1/2}}.
\end{eqnarray*}
For $N_2=O(1)$ we have the bound by
\begin{eqnarray*}
\lesssim\sum_{N_1,N_2=O(1)}\|P_{N_1}u\|_{L_{xT}^4}\|P_{N_2}u\|_{L_{xT}^4}\lesssim T^{\frac12}\|u\|_{X_T^{1/2}}^2,
\end{eqnarray*}
and the claim is proved.

It will thus suffice to show
\begin{eqnarray*}
\left\|\sum_{N_2\gg 1}(P_{1\ll \cdot\ll N_2}uP_{N_2}u_x)\right\|_{L_{xT}^2}^2\lesssim \|P_{\ge 1}u_0\|_{H^{1/2}}^4+(\|u\|_{X_T^{1/2}}^2+\|u\|_{X_T^{1/2}}^4)\|P_{\ge 1}u\|_{X_T^{1/2}}^2.
\end{eqnarray*}
(we take the square of $\|(u^2)_x\|_{L_{xT}^2}$.)
From the Littlewood-Paley theorem, we deduce the estimate
\begin{eqnarray*}
\left\|\sum_{N_2\gg 1}(P_{1\ll\cdot\ll N_2}uP_{N_2}u_x)\right\|_{L_{xT}^2}^2\sim \sum_{N_2\gg 1}\|P_{1\ll\cdot\ll N_2}uP_{N_2}u_x\|_{L_{xT}^2}^2
\end{eqnarray*}
which is written as
\begin{eqnarray}\label{eq:quard}
=\int_0^T\int_{-\infty}^{\infty}\sum_{1\ll N_1,N_1^*\ll N_2}P_{N_1}uP_{N_1^*}uP_{N_2}u_xP_{N_2}u_x\,dxdt.
\end{eqnarray}
We split the sum in $N_1,N_1^*$ as $\sum_{N_1\sim N_1^*}+\sum_{N_1\not\sim N_1^*}$.
The treatment for $N_1\sim N_1^*$ is as follows:
\begin{eqnarray*}
\lesssim\sum_{1\ll N_1\sim N_1^*\ll N_2}\|P_{N_1}u\|_{L_x^2L_T^{\infty}}\|P_{N_1^*}u\|_{L_x^2L_T^{\infty}}\|P_{N_2}u_x\|_{L_x^{\infty}L_T^2}^2,
\end{eqnarray*}
which is acceptable.
In order to study the contribution of $N_1\not\sim N_1^*$ for (\ref{eq:quard}),  we use the equation (\ref{eq:gBO}) to see that
\begin{eqnarray*}
(e^{it|\xi|\xi}\widehat{P_Nu}(\xi))_t=-e^{it|\xi|\xi}\widehat{P_N(u^2u_x)}(\xi).
\end{eqnarray*}
Then by Plancherel' theorem we can reduce to
\begin{eqnarray*}
=& & \sum_{1\ll N_1^*\ll N_1\ll N_2}\int_0^T\int_*e^{-it(|\xi_1|\xi_1+|\xi_2|\xi_2\pm\xi_3^2\mp\xi_4^2)}\\
& & (e^{it|\xi_1|\xi_1}\widehat{P_{N_1}u}(\xi_1))(e^{it|\xi_2|\xi_2}\widehat{P_{N_1^*}u}(\xi_2))(e^{\pm it\xi_3^2}\widehat{P_{\pm}P_{N_2}u_x}(\xi_3))(e^{\mp it\xi_4^2}\widehat{P_{\mp}P_{N_2}u_x}(\xi_4)),
\end{eqnarray*}
where we denote by $\int_*$ the integral over the hyper plane $\xi_{1234}=0$ (by symmetry we take $N_1^*\ll N_1$).
(Note that under the restriction $\xi_{1234}=0$, the terms of equal signs on $\xi_3,\xi_4$ vanish for $N_1^*\ll N_1\ll N_2$.)
Integrating by parts, we write this as
\begin{eqnarray}
&  & i\sum_{1\ll N_1^*\ll N_1\ll N_2}\int_*\frac{[\widehat{P_{N_1}u}(\xi_1,t)\widehat{P_{N_1^*}u}(\xi_2,t)\widehat{P_{\pm}P_{N_2}u_x}(\xi_3,t)\widehat{P_{\mp}P_{N_2}u_x}(\xi_4,t)]_{t=0}^{t=T}}{|\xi_1|\xi_1+|\xi_2|\xi_2\pm\xi_3^2\mp\xi_4^2}\label{eq:quard-1}\\
& & -i\sum_{1\ll N_1^*\ll N_1\ll N_2}\int_0^T\int_*\frac{1}{|\xi_1|\xi_1+|\xi_2|\xi_2\pm\xi_3^2\mp\xi_4^2}
\label{eq:quard-2}
\\
& & \left(\widehat{P_{N_1}(u^2u_x)}(\xi_1)\widehat{P_{N_1^*}u}(\xi_2)+\widehat{P_{N_1}u}(\xi_1)\widehat{P_{N_1^*}(u^2u_x)}(\xi_2)\right)\widehat{P_{\pm}P_{N_2}u_x}(\xi_3)\widehat{P_{\mp}P_{N_2}u_x}(\xi_4)
\nonumber
\\
& & -i\sum_{1\ll N_1^*\ll N_1\ll N_2}\int_0^T\int_*\frac{1}{|\xi_1|\xi_1+|\xi_2|\xi_2\pm\xi_3^2\mp\xi_4^2}\widehat{P_{N_1}u}(\xi_1)\widehat{P_{N_1^*}u}(\xi_2)
\label{eq:quard-3}
\\
& & \left(\widehat{P_{\pm}P_{N_2}(u^2u_x)_x}(\xi_3)\widehat{P_{\mp}P_{N_2}u_x}(\xi_4)+\widehat{P_{\pm}P_{N_2}u_x}(\xi_3)\widehat{P_{\mp}P_{N_2}(u^2u_x)_x}(\xi_4)\right),
\nonumber
\end{eqnarray}
where we omit the time variable $t$ for the sake of simplicity.

Observe that since $|\xi_2|\ll |\xi_1|\sim|\xi_{12}|\ll |\xi_3|\sim|\xi_4|$
\begin{eqnarray}\label{eq:deno}
\left|\frac{1}{|\xi_1|\xi_1+|\xi_2|\xi_2\pm\xi_3^2\mp\xi_4^2}+\frac{1}{\pm 2\xi_{12}\xi_3}\right|
\lesssim\frac{1}{\xi_3^2}.
\end{eqnarray}
Then by Coifman-Meyer's multilinear theorem \cite{CMS} (also in \cite{S}), the first term of the above integral (\ref{eq:quard-1}) is bounded by
\begin{eqnarray*}
& & \lesssim \sum_{t=0,T}\sum_{1\ll N_1^*\ll N_1\ll N_2}\\
& & (\|P_{N_1}u(t)\|_{L_x^4}\|P_{N_1^*}u(t)\|_{L_x^4}\|P_{N_2}u(t)\|_{L_x^4}^2+\|D_x^{-1}(P_{N_1}uP_{N_1^*}u)(t)\|_{L_x^{\infty}}\|D_x^{\frac12}P_{N_2}u(t)\|_{L_x^2}^2)\\
& & \lesssim \|P_{\ge 1}u_0\|_{H^{1/2}}^4+\|P_{\ge 1}u(T)\|_{H^{1/2}}^4
\end{eqnarray*}
which is acceptable.

The second term (\ref{eq:quard-2}) is treated in the same way as above.
We bound this contribution by
\begin{eqnarray*}
\lesssim\sum_{1\ll N_1^*\ll N_1\ll N_2}(\|P_{N_1}(u^3)\|_{L_{xT}^4}\|P_{N_1^*}u\|_{L_{xT}^4}+\|P_{N_1}u\|_{L_{xT}^4}\|P_{N_1^*}(u^3)\|_{L_{xT}^4})\|D_x^{\frac12}P_{N_2}u\|_{L_{xT}^4}^2
\end{eqnarray*}
which is easily acceptable.

To estimate the last term (\ref{eq:quard-3}), by (\ref{eq:deno}) we may replace the denominator in the integral term by $\pm 2i\xi_{12}\xi_3$.
This is because that if the denominator was $\xi_3^2$, we would have a bound by
\begin{eqnarray*}
\lesssim\sum_{1\ll N_1^*\ll N_1\ll N_2}\|P_{N_1}u\|_{L_{xT}^8}\|P_{N_1^*}u\|_{L_{xT}^8}\|D_x^{\frac12}P_{N_2}u\|_{L_{xT}^4}\|D_x^{\frac12}P_{N_2}(u^3)\|_{L_{xT}^2}.
\end{eqnarray*}
To estimate the last term $\|D_x^{\frac12}P_{N_2}(u^3)\|_{L_{xT}^2}$, we use the Cauchy-Schwartz inequality in $N_2$ to handle this by using that
\begin{eqnarray*}
\left(\sum_{N_2}\|D_x^{1/2}P_{N_2}(u^3)\|_{L_{xT}^2}^2\right)^{1/2}\lesssim\|D_x^{\frac12}(u^3)\|_{L_{xT}^2}\lesssim \|D_x^{\frac12}u\|_{L_{xT}^4}\|u\|_{L_{xT}^8}^2\lesssim \|u\|_{X_T^{1/2}}^3.
\end{eqnarray*}
We can use symmetry in (\ref{eq:quard-3}) because $-\xi_3=\xi_4+\xi_{12}$, and thus we are reduced to the following integral
\begin{eqnarray*}
\sum_{1\ll N_1^*\ll N_1\ll N_2}\int_0^T\int_{*}\left(2\frac{\xi_3\xi_4}{\xi_{12}}+\xi_3\right)\widehat{P_{N_1}u}(\xi_1)\widehat{P_{N_1^*}u}(\xi_2)\widehat{P_{\pm}P_{N_2}(u^3)}(\xi_3)\widehat{P_{\mp}P_{N_2}u}(\xi_4).
\end{eqnarray*}
For the Fourier multiplier $\xi_3$, by H\"older's inequality we obtain
\begin{eqnarray*}
\lesssim\sum_{1\ll N_1^*\ll N_1\ll N_2}\|P_{N_1}u\|_{L_{xT}^8}\|P_{N_1^*}u\|_{L_{xT}^8}\|D_x^{\frac12}P_{N_2}u\|_{L_{xT}^4}\|D_x^{\frac12}P_{N_2}(u^3)\|_{L_{xT}^2},
\end{eqnarray*}
which can be treated as before.
On the other hand, for the Fourier multiplier $2\frac{\xi_3\xi_4}{\xi_{12}}$, we deduce from H\"older's inequality that the expression is bounded by
\begin{eqnarray*}
\lesssim\sum_{1\ll N_1^*\ll N_1\ll N_2}\|D_x^{-1}(P_{N_1}uP_{N_1^*}u)\|_{L_x^2L_T^{\infty}}\|\partial_xP_{N_2}u\|_{L_x^{\infty}L_T^2}\|\partial_xP_{N_2}(u^3)\|_{L_{xT}^2}.
\end{eqnarray*}
(We may of course decompose $P_{N_1}u=P_+P_{N_1}u+P_-P_{N_1}u$ in order to obtain the $D_x^{-1}$-derivative.)
By Young's inequality, we see that the first term is bounded by
\begin{eqnarray*}
& & \sum_{1\ll N_1^*\ll N_1}\|D_x^{-1}(P_{N_1}uP_{N_1^*}u)\|_{L_x^2L_T^{\infty}}\\
&\lesssim & \sum_{1\ll N_1^*\ll N_1}\|[D_x^{-1}\widetilde{P_{N_1}}]\check{~}\|_{L_x^2}\|P_{N_1}uP_{N_1^*}u\|_{L_x^1L_T^{\infty}}\\
&\lesssim & \sum_{1\ll N_1^*\ll N_1}\|\langle D_x\rangle^{-\frac14}P_{N_1}u\|_{L_x^2L_T^{\infty}}\|\langle D_x\rangle^{-\frac14}P_{N_1^*}u\|_{L_x^2L_T^{\infty}}\\
& \lesssim &\|P_{\ge 1}u\|_{X_T^{1/2}}^2.
\end{eqnarray*}
The second term $\|\partial_xP_{N_2}u\|_{L_x^{\infty}L_T^2}$ yields an acceptable term after Cauchy-Schwarz.
Also for the last term $\|\partial_xP_{N_2}(u^3)\|_{L_{xT}^2}$, after Cauchy-Schwartz, we control this by
\begin{eqnarray*}
\lesssim\left(\sum_{N_2}\|\partial_xP_{N_2}(u^3)\|_{L_{xT}^2}^2\right)^{\frac12}\sim\|(u^3)_x\|_{L_{xT}^2}.
\end{eqnarray*}
We use a Littlewood-Paley decomposition to expand
\begin{eqnarray*}
u^3 & = & \sum_{N\ge 1}\left((P_{\le 2N}u)^3-(P_{\le N}u)^3\right)+P_0(u)^3\\
& =& \sum_{N\ge 1}(P_{\le 2N}-P_{\le N})uA_N(u)+P_0(u)^3\\
& =& \sum_{N\ge 1}P_{2N}uA_N(u)+P_0(u)^3,
\end{eqnarray*}
where $A_N(u)=(P_{\le 2N}u)^2+P_{\le 2N}uP_{\le N}u+(P_{\le N}u)^2$.
The term $P_0(u)^3$ can be estimated by using a Littlewood-Paley theorem to estimate
\begin{eqnarray*}
\|\partial_xP_0(u)^3\|_{L_{xT}^2}\lesssim\|P_0(u)\|_{L_{xT}^6}^3\lesssim\|u\|_{X_T^{1/2}}^3.
\end{eqnarray*}
For the first term, suppose $A_N(u)=(P_{\le N}u)^2$, because other two terms are handled similarly.
Then it suffices to show
\begin{eqnarray*}
\left\|\partial_x\sum_{N\ge 1}P_{2N}u(P_{\le N}u)^2\right\|_{L_{xT}^2}\lesssim \|u\|_{X_T^{1/2}}^3.
\end{eqnarray*}
Now we write
\begin{eqnarray*}
\partial_x\sum_{N\ge 1}P_{2N}u(P_{\le N}u)^2
=\left(\sum_{M\lesssim 1}+\sum_{M\gg 1}\right)\partial_xP_M\sum_{N\ge 1}P_{2N}u(P_{\le N}u)^2.
\end{eqnarray*}

The low frequency part $\sum_{M\lesssim 1}$ can be estimated by using a Littlewood-Paley theorem to estimate
\begin{eqnarray*}
\left\|\partial_xP_{\lesssim 1}\sum_{N\ge 1}P_{2N}u(P_{\le N}u)^2\right\|_{L_{xT}^2}& \lesssim &\left\|P_{\lesssim 1}\sum_{N\ge 1}P_{2N}u(P_{\le N}u)^2\right\|_{L_{xT}^2}\\
& \lesssim & \sum_{N\ge 1}\|P_{2N}u(P_{\le N}u)^2\|_{L_{xT}^2}\lesssim \sum_{N\ge 1}\|P_{2N}u\|_{L_{xT}^6}\|u\|_{L_{xT}^6}^2.
\end{eqnarray*}
But by
\begin{eqnarray*}
\sum_{N\ge 1}\|P_{2N}u\|_{L_{xT}^6}\le\left(\sum_{N\ge 1}\|P_Nu\|_{L_T^6W_x^{1/2,6}}^2\right)^{1/2}\lesssim \|u\|_{X_T^{1/2}},
\end{eqnarray*}
we can bound the left-hand side by $\lesssim \|u\|_{X_T^{1/2}}^3$.

We now look at the contribution of the sum $\sum_{M\gg 1}$.
We begin by using a Littlewood-Paley theorem to write
\begin{eqnarray*}
\left\|\sum_{M\gg 1}\partial_xP_M\sum_{N\ge 1}P_{2N}u(P_{\le N}u)^2\right\|_{L_{xT}^2} & \sim & \left(\sum_{M\gg 1}M^2\|P_M\sum_{N\ge 1}P_{2N}u(P_{\le N}u)^2\|_{L_{xT}^2}^2\right)^{1/2}\\
& \lesssim & \left(\sum_{M\gg 1}M^2\|P_M\sum_{N\gtrsim M}P_{2N}u(P_{\le N}u)^2\|_{L_{xT}^2}^2\right)^{1/2}.
\end{eqnarray*}
Using H\"older and Littlewood-Paley theorem, this is bounded by
\begin{eqnarray*}
& & \lesssim \left(\sum_{M\gg 1}\left(\sum_{N\gtrsim M}\frac{M}{N}\|P_{2N}u_x\|_{L_x^{\infty}L_T^2}\|P_{\le N}u\|_{L_x^4L_T^{\infty}}^2\right)^2\right)^{1/2}.
\end{eqnarray*}
Now we have
\begin{eqnarray*}
\sup_N\|P_{\le N}u\|_{L_x^4L_T^{\infty}}\lesssim\sup_N\sum_{M\le N}\|P_Mu\|_{L_x^4L_T^{\infty}}\lesssim \|u\|_{X_T^{1/2}}
\end{eqnarray*}
and hence
\begin{eqnarray*}
\left(\sum_{M\gg 1}\left(\sum_{N\gtrsim M}\frac{M}{N}\|P_{2N}u_x\|_{L_x^{\infty}L_T^2}\|P_{\le N}u\|_{L_x^4L_T^{\infty}}^2\right)^2\right)^{1/2}
\lesssim\|u\|_{X_T^{1/2}}^2\left(\sum_{M\gg 1}\left(\sum_{N\gtrsim N}\frac{M}{N}\|P_{N}u_x\|_{L_x^{\infty}L_T^2}\right)^2\right)^{1/2}.
\end{eqnarray*}
But using Young's inequality $\|f*g\|_{l^2}\le\|f\|_{l^1}\|g\|_{l^2}$, we bound the left-hand side by
\begin{eqnarray*}
\lesssim \|u\|_{X_{T}^{1/2}}^2\left(\sum_{N\gg 1}\|P_Nu_x\|_{L_x^{\infty}L_T^2}^2\right)^{1/2}\le\|u\|_{X_T^{1/2}}^3. 
\end{eqnarray*}
This concludes the proof of Proposition \ref{prop:L^2-spacetime}.
\qed 
\section{Nonlinear estimates}\label{sec:nonlinear}
We shall now deal with the problem of estimating the nonlinearity arising in the equation (\ref{eq:mmBO5}).
Throughout this section, we always assume $T\in (0,1)$ and $N\gg 1$.
For brevity's sake, we only consider the endpoint case $s=1/2$, and abbreviate $X_T^{1/2}$ to $X$.
Recall that the equation (\ref{eq:mmBO5}) has the following equivalent integral equation
\begin{eqnarray*}
v_N(t) & = & e^{it\partial_x^2}(e^{-\frac{i}{2}\int^x(P_{\ll N}u_0)^2}P_+P_Nu_0)\label{eq:I-BO1}\\
& & +\int_0^te^{i(t-t')\partial_x^2}(A_{1,N}+A_{2,N}+A_{3,N}+A_{4,N}+A_{5,N})(t')\,dt',
\end{eqnarray*}
where $A_{j,N}(t)$ are defined in (\ref{eq:mmBO5}).

Because of Lemmas \ref{lem:freeestimate} and \ref{lem:retard}, we need to define the function space $Y$, equipped with the following norm, which will only be used in this section and next section (Lemmas \ref{lem:low-u} and \ref{lem:u-v}) 
\begin{eqnarray*}
\|u\|_Y=\|u\|_{L_T^{\infty}H_x^{1/2}}+\|\partial_xu\|_{L_x^{\infty}L_T^2}+\|u\|_{L_x^2L_T^{\infty}}+\|\langle D_x\rangle^{1/4}u\|_{L_x^4L_T^{\infty}}.
\end{eqnarray*}
We handle the $Y$-norm for the nonlinearities $A_j$.
\begin{proposition}\label{prop:I-BO7}
Let $u$ be a $H^{\infty}$-solution to (\ref{eq:gBO}).
Then
\begin{eqnarray*}
& & \left(\sum_{N\gg 1}\left\|\int_0^te^{i(t-t')\partial_x^2}\sum_{j=1}^5A_{j,N}(t')\,dt'
\right\|_Y^2\right)^{1/2}
\nonumber\\
&  & \lesssim T^{1/4}(1+\|u\|_X^4)\|P_{\gtrsim 1}u\|_X+\|P_{\gtrsim 1}u_0\|_{H^{1/2}}^2\|P_{\gtrsim 1}u\|_X+(1+\|u\|_X^2)\|P_{\gtrsim 1}u\|_X^2.
\label{eq:nI-BO2}
\end{eqnarray*}
\end{proposition}
{\it Proof of Proposition \ref{prop:I-BO7}.}
We consider each contribution separately.
\subsection{The contribution of $A_{1,N}$.}
We begin with the identity
\begin{eqnarray*}
u^2u_x=(P_{\ll N}u)^2\widetilde{P}_Nu_x+(P_{\ll N}u)^2(1-\widetilde{P}_N)u_x+(u^2-(P_{\ll N}u)^2)u_x.
\end{eqnarray*}
Note that the term $(P_{\ll N}u)^2P_+\widetilde{P}_Nu_x$ has Fourier support in $|\xi|\sim N$, also the second term $(P_{\ll N}u)^2(1-\widetilde{P}_N)u_x$ will cancel when the projection operator $P_N$ is applied since $P_N((P_{\ll N}u)^2u_x)=P_N((P_{\ll N}u)^2\widetilde{P}_Nu_x)$ for large frequency $\sim N$.
On the other hand, $P_+((P_{\ll N}u)^2\widetilde{P}_Nu_x)=(P_{\ll N}u)^2P_+\widetilde{P}_Nu_x$ for large frequency $\sim N$.
We thus have
\begin{eqnarray}
& & P_+P_N(u^2u_x)-(P_{\ll N}u)^2P_+P_Nu_x
\nonumber\\
& = & P_N((P_{\ll N}u)^2P_+\widetilde{P}_Nu_x)-(P_{\ll N}u)^2P_NP_+\widetilde{P}_Nu_x
\label{eq:D1}\\
& & +P_+P_N((u^2-(P_{\ll N}u)^2)u_x),
\label{eq:D2}
\end{eqnarray}
where we may freely add $\widetilde{P}_N$ to $P_N$.
We exploit the projection operator $P_N$ to expand the second term (\ref{eq:D2}) as follows: for each $N\gg 1$
\begin{eqnarray*}
P_+P_N((u^2-(P_{\ll N}u)^2)u_x)=P_+P_N\sum_{k=1}^5\sum_{I_k}P_{N_1}uP_{N_2}uP_{N_3}u_x,
\end{eqnarray*}
where
\begin{eqnarray}
& & I_1:~N_1\sim N_2\gtrsim N,N_3,\label{eq:I_1}\\
& & I_2:~N_1\sim N\gg N_2,N_3,\label{eq:I_2}\\
& & I_3:~N_1\sim N_3\sim N\gg N_2,\label{eq:I_3}\\
& & I_4:~N_1\sim N_3\gg N_2\gtrsim N,\label{eq:I_4}\\
& & I_5:~N_1\sim N_3\gg N\gg N_2.\label{eq:I_5}
\end{eqnarray}
(By symmetry, we may assume $N_1\gtrsim N_2,N_3,N$.)

We now give the following lemma.
\begin{lemma}\label{lem:L^2}
Let $u$ be a solution of (\ref{eq:gBO}).
Then
\begin{eqnarray}
& & \left(\sum_{N\gg 1}\|P_N((P_{\ll N}u)^2P_+\widetilde{P}_Nu_x)-(P_{\ll N}u)^2P_NP_+\widetilde{P}_Nu_x\|_{L_x^1L_T^2}^2\right)^{1/2}
\label{eq:L^2-1}\\
& & +\left(\sum_{N\gg 1}\sum_{k=2,3,5}\|P_N\sum_{I_k}P_{N_1}uP_{N_2}uP_{N_3}u_x\|_{L_x^1L_T^2}^2\right)^{1/2}
\label{eq:L^2-2}\\
& \lesssim & T^{1/4}\|u\|_X^2\|P_{\gtrsim 1}u\|_X+\|P_{\gtrsim 1}u_0\|_{H^{1/2}}^2\|P_{\gtrsim 1}u\|_{X}+(1+\|u\|_X^2)\|P_{\gtrsim 1}u\|_X^2.
\nonumber
\end{eqnarray}
\end{lemma}
{\it Proof of Lemma \ref{lem:L^2}.}
We  first consider (\ref{eq:L^2-1}).
To shift a derivative from the high-frequency function $P_+\widetilde{P}_Nu_x$ to the low-frequency function $(P_{\ll N}u)^2$, we require the following Leibniz rule for $P_N$
\begin{eqnarray*}\label{eq:P_Leibniz}
\widehat{(P_N(fg)-fP_Ng)}(\xi)&=&i\int(\varphi_N(\xi)-\varphi_N(\xi_1))\hat{f}(\xi-\xi_1)\hat{g}(\xi_1)\,d\xi_1\\
&=& i\int\left(\int_0^1\varphi_N'((1-\eta)\xi_1+\eta\xi)\,d\eta\right)(\xi-\xi_1)\hat{f}(\xi-\xi_1)\hat{g}(\xi_1)\,d\xi_1,
\end{eqnarray*}
and its Fourier inverse formula
\begin{eqnarray*}\label{eq:I_P}
(P_N(fg)-fP_Ng)(x)=\int_0^1\,d\eta\left(\int\check{\varphi}_N(y)yf_x(x-\eta y)g(x-y)\,dy\right).
\end{eqnarray*}
Since $\|y\check{\varphi}_N\|_{L_y^1}=cN^{-1}\|y\check{\varphi}_1\|_{L_y^1}$, we may bound the contribution of (\ref{eq:L^2-1}) by
\begin{eqnarray*}
& \lesssim & (\sum_{N\gg 1}\|(P_{\ll N}u)^2_x\|_{L_{xT}^2}^2N^{-2}\|P_+P_Nu_x\|_{L_x^{2}L_T^{\infty}}^2)^{1/2}\\
& \lesssim & (\sum_{N\gg 1}\|(P_{\ll N}u)^2_x\|_{L_{xT}^2}^2\|P_Nu\|_{L_x^{2}L_T^{\infty}}^2)^{1/2}.
\end{eqnarray*}
Split $P_{\ll N}u=P_{\le 1}u+P_{1<\cdot\ll N}u$, and write $(P_{\ll N}u)^2=(P_{\le 1}u)^2+2P_{\le 1}uP_{1<\cdot\ll N}u+(P_{1<\cdot\ll N}u)^2$.
For $(P_{\le 1}u)^2$, we can discard the $\partial_x$-derivative, and estimate this contribution by
\begin{eqnarray*}
\lesssim T^{\frac12}\|u\|_X^2\|P_{\gtrsim 1}u\|_X.
\end{eqnarray*}
For the contributions of the other two terms, we use Proposition \ref{prop:L^2-spacetime}\footnote{More precisely, we use the proof of Proposition \ref{prop:L^2-spacetime} and replace $(u^2)_x$ with $(P_{\le 1}uP_{1<\cdot\ll N}u)_x$ or $((P_{1<\cdot\ll N}u)^2)_x$.} to obtain the desired bound.

Turning to the estimate (\ref{eq:L^2-2}), we shall consider separately the contributions of (\ref{eq:I_2}), (\ref{eq:I_3}) and (\ref{eq:I_5}).
For (\ref{eq:I_2}), we bound this contribution  to (\ref{eq:L^2-2}) by
\begin{eqnarray*}
& & \lesssim\left(\sum_{N\gg 1}\|\sum_{N_1\sim N}P_{N_1}uP_{\ll N}uP_{\ll N}u_x\|_{L_x^1L_T^2}^2\right)^{1/2}\\
& & \lesssim\left(\sum_{N\gg 1}\|P_Nu\|_{L_x^2L_T^{\infty}}^2\|P_{\ll N}uP_{\ll N}u_x\|_{L_{xT}^2}^2\right)^{1/2},
\end{eqnarray*}
which is acceptable, since the proof is along the same lines as that for (\ref{eq:L^2-1}).
For (\ref{eq:I_3}), a similar argument shows that this contribution to (\ref{eq:L^2-2}) is bounded by
\begin{eqnarray*}
& \lesssim & \left(\sum_{N\gg 1}\|\sum_{N_1\sim N}(P_{N_1}u\sum_{N_3\sim N_1}P_{\ll N}uP_{N_3}u_x)\|_{L_x^1L_T^2}^2\right)^{1/2}\\
& \lesssim & \left(\sum_{N\gg 1}\left(\sum_{N_1\sim N}\|P_{N_1}u\sum_{N_3\sim N_1}P_{\ll N}uP_{N_3}u_x\|_{L_x^1L_T^2}\right)^2\right)^{1/2}\\
& \lesssim & \left(\sum_{N\gg 1}\sum_{N_1\sim N}\|P_{N_1}u\|_{L_x^2L_T^{\infty}}^2\sum_{N_1\sim N}\|\sum_{N_3\sim N_1}P_{\ll N}uP_{N_3}u_x\|_{L_{xT}^2}^2\right)^{1/2}.
\end{eqnarray*}
This is bounded by
\begin{eqnarray*}
\lesssim\left(\sum_{N\gg 1}\|P_{N}u\|_{L_x^2L_T^{\infty}}^2\|P_{\ll N}uP_{\sim N}u_x\|_{L_{xT}^2}^2\right)^{1/2},
\end{eqnarray*}
which is acceptable as before.
Finally, for (\ref{eq:I_5}), we observe that by symmetry
\begin{eqnarray*}
\sum_{N_1\sim N_3}P_{N_1}uP_{N_3}u_x=\sum_{N_1\sim N_3}(P_{N_1}uP_{N_3}u)_x.
\end{eqnarray*}
Then we bound this contribution to (\ref{eq:L^2-2}) by
\begin{eqnarray*}
& & \lesssim\left(\sum_{N\gg 1}\|\widetilde{P}_N(P_{\gg N}uP_{\gg N})_xP_{\ll N}u\|_{L_x^1L_T^2}^2\right)^{1/2}\\
& & \lesssim\left(\sum_{N\gg 1}N^{2}\sum_{N_1\gg N}\|P_{N_1}u\|_{L_x^2L_T^{\infty}}^2\sum_{N_1\gg N}\|\sum_{N_3\sim N_1}P_{N_3}uP_{\ll N}u\|_{L_{xT}^2}^2\right)^{1/2}\\
& & \lesssim\|P_{\gtrsim 1}u\|_{X}\left(\sum_{N\gg 1}N^{2}\sum_{N_3\gg N}N_3^{-2}\|P_{N_3}u_xP_{\ll N}u\|_{L_{xT}^2}^2\right)^{1/2}.
\end{eqnarray*}
We rearrange the sum as follows
\begin{eqnarray}
& \lesssim & \|P_{\gtrsim 1}u\|_{X}\left(\sum_{N_3\gg N\gg 1}(N/N_3)^{2}\|P_{N_3}u_xP_{\ll N}u\|_{L_{xT}^2}^2\right)^{1/2}\nonumber\\
& \lesssim & \|P_{\gtrsim 1}u\|_{X}\left(\sum_{N_3\gg 1}\sum_{N:N\ll N_3}(N/N_3)^{2}\sup_{N:N\ll N_3}\|P_{N_3}u_xP_{\ll N}u\|_{L_{xT}^2}^2\right)^{1/2}.\label{eq:mL^2}
\end{eqnarray}
One can then observe the following variant of (\ref{eq:L^2-spacetime}) that entered in the proof of Proposition \ref{prop:L^2-spacetime}\footnote{Incidentally, the estimate (\ref{eq:rrL^2}) holds without the last term $T^{1/2}\|u\|_X^2$ under the restriction $N_1\ll N_2$}:
\begin{eqnarray}
& & \left(\sum_{N_2\gg 1}\sup_{N_1:N_1\ll N_2}\|P_{N_2}u_xP_{\ll N_1}u\|_{L_{xT}^2}^2\right)^{1/2}
\nonumber\\
&\lesssim & \|P_{\gg 1}u_0\|_{H^{1/2}}^2+(1+\|u\|_X)\|u\|_X\|P_{\gg 1}u\|_X+T^{1/2}\|u\|_X^2.
\label{eq:rrL^2}
\end{eqnarray}
The idea is that the contribution of the term $P_{\ll N_1}u$ can be essentially estimated by the squared-type norm $(\sum_{M\ll N_1}\|P_{M}u\|_X^2)^{1/2}\lesssim \|u\|_{X}$, which is independent of the size of $N_1$.
We first sum in $N$, then in $N_3$ for (\ref{eq:mL^2}), and use the inequality (\ref{eq:rrL^2}).

This completes the proof of Lemma \ref{lem:L^2}.
\qed

We now turn to the proof of Proposition \ref{prop:I-BO7}, and estimate the contribution of $A_{1,N}$.
With the aid of this lemma, we can prove the estimate for the terms (\ref{eq:D1}), (\ref{eq:I_1}), (\ref{eq:I_2}), (\ref{eq:I_3}), (\ref{eq:I_4}) and (\ref{eq:I_5}).
We shall consider separately each contribution.

\subsubsection{The contribution of (\ref{eq:D1}).}
By (\ref{eq:rstrichartz-2}), (\ref{eq:rsmoothing-2}), (\ref{eq:rmaximal-1}), (\ref{eq:rmaximal-2}), we bound the contribution of (\ref{eq:D1}) to the left of (\ref{eq:nI-BO2}) by
\begin{eqnarray}
&  & \lesssim \left(\sum_{N\gg 1}\|B_N(u)\|_{L_x^1L_T^2}^2\right)^{1/2}
\label{eq:D1-1}\\
& & +\left(\sum_{N\gg 1}\left(\sum_M\|P_M(e^{-\frac{i}{2}\int^x(P_{\ll N}u)^2}B_N(u))\|_{L_x^1L_T^2}\right)^2\right)^{1/2},
\label{eq:D1-2}
\end{eqnarray}
where $B_N(u)=P_N((P_{\ll N}u)^2P_+\widetilde{P}_Nu_x)-(P_{\ll N}u)^2P_NP_+\widetilde{P}_Nu_x$.
From Lemma \ref{lem:L^2}, the first term (\ref{eq:D1-1}) is acceptable.

On the other hand, for the second term (\ref{eq:D1-2}), we split the sum $\sum_M$ into three parts $\sum_{M\sim N}+\sum_{M\ll N}+\sum_{M\gg N}$.
The contribution of $M\sim N$ is of type (\ref{eq:D1-1}) by summing in $M$ such that $M\sim N$.
Next we study the contribution of $M\ll N$ to (\ref{eq:D1-2}).
Since the expression $P_N((P_{\ll N}u)^2P_+\widetilde{P}_Nu_x)-(P_{\ll N}u)^2P_NP_+\widetilde{P}_Nu_x$ has Fourier support in $|\xi|\sim N$, we may add the projection operator $P_{\sim N}$ to $e^{-\frac{i}{2}\int_{-\infty}^x(P_{\ll N}u)^2\,dy}$.
By H\"older inequality, we can bound this contribution to (\ref{eq:D1-2}) by
\begin{eqnarray*}
& & \lesssim \left(\sum_{N\gg 1}\left(\sum_{M\ll N}\|P_{\sim N}e^{-\frac{i}{2}\int^x(P_{\ll N}u)^2}B_N(u)\|_{L_x^1L_T^2}\right)^2\right)^{1/2}\\
& & \lesssim \left(\sum_{N\gg 1}(\log N)^2\|P_{\sim N}e^{-\frac{i}{2}\int^x(P_{\ll N}u)^2}\|_{L_{xT}^{1/\varepsilon}}^2\|B_N(u)\|_{L_x^{\frac{1}{1-\varepsilon}}L_T^{\frac{2}{1-2\varepsilon}}}^2\right)^{1/2}.
\end{eqnarray*}
We easily see that by Sobolev inequality
\begin{eqnarray*}
N\|P_{\sim N}e^{-\frac{i}{2}\int^x(P_{\ll N}u)^2}\|_{L_{xT}^{1/\varepsilon}}& \lesssim & \|\partial_xP_{\sim N}e^{-\frac{i}{2}\int^x(P_{\ll N}u)^2\,dy}\|_{L_{xT}^{1/\varepsilon}}\\
&\lesssim&\|P_{\ll N}u\|_{L_{xT}^{\frac{2}{\varepsilon}}}^2\lesssim\|u\|_X^2,
\end{eqnarray*}
and
\begin{eqnarray*}
\sum_{N\gg 1}N^{-2+2\varepsilon}\|B_N(u)\|_{L_x^{\frac{1}{1-\varepsilon}}L_T^{\frac{2}{1-2\varepsilon}}}^2
&\lesssim & \sum_{N\gg 1}\|P_{\ll N}u\|_{L_x^{4}L_T^{\infty}}^4\|D_x^{\varepsilon}P_Nu\|_{L_{xT}^{\frac{2}{1-2\varepsilon}}}^2\\
&\lesssim & T^{1-2\varepsilon}\|u\|_X^4\|P_{\gtrsim 1}u\|_X^2.
\end{eqnarray*}
From these, the previous is bounded by
\begin{eqnarray*}
\lesssim T^{\frac12-\varepsilon}\|u\|_X^4\|P_{\gtrsim 1}u\|_X.
\end{eqnarray*}

For the contribution of $M\gg N$ to (\ref{eq:D1-2}), we now add the projection operator $P_M$ to $e^{-\frac{i}{2}\int_{-\infty}^x(P_{\ll N}u)^2\,dy}$.
Then we have the bound by
\begin{eqnarray*}
\lesssim \left(\sum_{N\gg 1}\left(\sum_{M\gg N}\|P_{M}e^{-\frac{i}{2}\int^x(P_{\ll N}u)^2}\|_{L_{xT}^{1/\varepsilon}}\|B_N(u)\|_{L_x^{\frac{1}{1-\varepsilon}}L_T^{\frac{2}{1-2\varepsilon}}}\right)^2\right)^{1/2}
\end{eqnarray*}
and this follows from the same line of proof as the contribution of the case $M\ll N$.
This completes the proof for (\ref{eq:D1}).

\subsubsection{The contribution of (\ref{eq:I_1}).}
We use (\ref{eq:rstrichartz-1}), (\ref{eq:srsmoothing}), (\ref{eq:srmaximal-1}), (\ref{eq:srmaximal-2}), and estimate this by
\begin{eqnarray}
& \lesssim & \left(\sum_{N\gg 1}\|e^{-\frac{i}{2}\int^x(P_{\ll N}u)^2}C_N(u)\|_{L_T^1H_x^{\frac12}}^2\right)^{1/2}
\label{eq:D2-1-1}\\
& & +\left(\sum_{N\gg 1}\left(\sum_M\|P_M(e^{-\frac{i}{2}\int^x(P_{\ll N}u)^2}C_N(u))\|_{L_T^1H_x^{\frac12}}\right)^2\right)^{1/2},
\label{eq:D2-1-2}
\end{eqnarray}
where $C_N(u)=P_+P_N\sum_{I_1}P_{N_1}uP_{N_2}uP_{N_3}u_x$.
By Lemma \ref{lem:Leibnitz2}, the first term (\ref{eq:D2-1-1}) is bounded by
\begin{eqnarray}
& \lesssim & T^{\frac12}\left(\sum_{N\gg 1}\|P_{\ll N}u\|_{L_{xT}^{\frac{4}{\varepsilon}}}^4\|P_N\sum_{I_1}P_{N_1}uP_{N_2}uP_{N_3}u_x\|_{L_{xT}^{\frac{2}{1-\varepsilon}}}^2\right)^{1/2}\label{eq:D2-1-3}\\
& & +T^{\frac12}\left(\sum_{N\gg 1}\|D_x^{\frac12}P_N\sum_{I_1}P_{N_1}uP_{N_2}uP_{N_3}u_x\|_{L_{xT}^2}^2\right)^{1/2}.\label{eq:D2-1-4}
\end{eqnarray}
It is easy to see that by Sobolev inequality the term (\ref{eq:D2-1-3}) is bounded by
\begin{eqnarray*}
\lesssim T^{1/2}\|u\|_X^4\|P_{\gtrsim 1}u\|_X.
\end{eqnarray*}
For the term (\ref{eq:D2-1-4}), we may drop the assumption on $N$ for (\ref{eq:I_1}), namely $I_1:N_1\sim N_2\gtrsim N_3$.
In fact, add
\begin{eqnarray*}
P_N\sum_{N\gg N_1\sim N_2\gtrsim N_3}P_{N_1}uP_{N_2}uP_{N_3}u_x=0.
\end{eqnarray*}
We therefore bound (\ref{eq:D2-1-4}) by
\begin{eqnarray*}
& & \lesssim T^{1/2}\|D_x^{1/2}\sum_{N_1\sim N_2\gtrsim N_3}P_{N_1}uP_{N_2}uP_{N_3}u_x\|_{L_{xT}^2}\\
& & \lesssim \sum_{N_1\sim N_2}\|\langle D_x\rangle^{1/2}P_{N_1}u\|_{L_{xT}^6}\|\langle D_x\rangle^{1/2}P_{N_2}u\|_{L_{xT}^6}\|\langle D_x\rangle^{1/2}P_{\lesssim N_1}u\|_{L_{xT}^6},
\end{eqnarray*}
so that summing on $N_1\sim N_2(\gtrsim N\gg 1)$ gives
\begin{eqnarray*}
\lesssim T^{1/2}\|u\|_X^2\|P_{\gtrsim 1}u\|_X.
\end{eqnarray*}
For (\ref{eq:D2-1-2}), we split the sum $\sum_M$ into two parts $\sum_{M\lesssim N}+\sum_{M\gg N}$, which gives the bound by
\begin{eqnarray*}
& \lesssim & \left(\sum_{N\gg 1}\left(\sum_{M\lesssim N}\langle M\rangle^{\frac12}\|C_N(u)\|_{L_T^1L_x^2}\right)^2\right)^{1/2}\\
& & +\left(\sum_{N\gg 1}\left(\sum_{M\gg N}\|P_M\langle D_x\rangle^{\frac12}e^{-\frac{i}{2}\int^x(P_{\ll N}u)^2}C_N(u)\|_{L_T^1L_x^2}\right)^2\right)^{1/2}.
\end{eqnarray*}
Since $\sum_{M\lesssim N}\langle M\rangle^{1/2}\lesssim N^{1/2}$, the estimate for the first term follows from the same argument as that for (\ref{eq:D2-1-4}).
The second term is treated by
\begin{eqnarray*}
\lesssim T^{\frac12}\left(\sum_{N\gg 1}\left(\sum_{M\gg N}M^{-\varepsilon}\right)^{1/2}\|\partial_xe^{-\frac{i}{2}\int^x(P_{\ll N}u)^2}\|_{L_{xT}^{\frac{2}{\varepsilon}}}^2\|C_N(u)\|_{L_{xT}^{\frac{2}{1-\varepsilon}}}^2\right)^{1/2}.
\end{eqnarray*}
By the same argument as in (\ref{eq:D2-1-3}), this is  bounded by
\begin{eqnarray*}
\lesssim T^{\frac12}\|u\|_X^4\|P_{\gtrsim 1}u\|_X.
\end{eqnarray*}
Then this gives the proof for the contribution of (\ref{eq:I_1}).

\subsubsection{The contribution of (\ref{eq:I_2}).}
It is useful to recall the proof for the contribution of (\ref{eq:D1}).
We use (\ref{eq:rstrichartz-2}), (\ref{eq:rsmoothing-2}), (\ref{eq:rmaximal-1}), (\ref{eq:rmaximal-2}) to obtain the bound by
\begin{eqnarray*}
&\lesssim & \left(\sum_{N\gg 1}\|\sum_{N_1\sim N}P_{N_1}uP_{\ll N}uP_{\ll N}u_x\|_{L_x^1L_T^2}^2\right)^{1/2}\\
& & +\left(\sum_{N\gg 1}\left(\sum_M\|P_M(e^{-\frac{i}{2}\int^x(P_{\ll N}u)^2}P_+P_N\sum_{N_1\sim N}P_{N_1}uP_{\ll N}uP_{\ll N}u_x)\|_{L_T^1L_x^2}\right)^2\right)^{1/2}.
\end{eqnarray*}
The proof for (\ref{eq:L^2-1}) in Lemma \ref{lem:L^2} leads to that for the first term.
The proof for the second term follows the corresponding argument for the contribution of (\ref{eq:D1}).

\subsubsection{The contribution of (\ref{eq:I_3}).}
This follows from the same argument as that for (\ref{eq:I_2}).

\subsubsection{The contribution of (\ref{eq:I_4}).}
We invoke the proof used for the contribution of (\ref{eq:I_1}) to prove the estimate.
We can share the derivative with three $P_{N_1}u,P_{N_2}u$ and $P_{N_3}u$, and also the derivative in $P_{N_3}u_x$ can be shifted to that for $P_{N_1}u$, in view of the support property $N_1\sim N_3\gg N_2\gtrsim N\gg 1$.
Hence the same proof as that for the case (\ref{eq:I_1}) gives us the desired conclusion.

\subsubsection{The contribution of (\ref{eq:I_5}).}
In order to verify the proof of (\ref{eq:I_5}), we reprise the proof of (\ref{eq:D1}), using Lemma \ref{lem:L^2}.
It thus remains only to estimate
\begin{eqnarray*}
\left(\sum_{N\gg 1}\left(\sum_M\|P_M(e^{-\frac{i}{2}\int^x(P_{\ll N}u)^2}P_+P_N\sum_{N_1\sim N_3\gg N}P_{N_1}uP_{\ll N}uP_{N_3}u_x)\|_{L_T^1L_x^2}\right)^2\right)^{1/2}.
\end{eqnarray*}
We again exploit the projection operator $P_N$ to obtain
\begin{eqnarray*}
& & \sum_{N\gg 1}N^{-2+2\varepsilon}\|\widetilde{P}_N(P_{\gg N}uP_{\gg N}u)_xP_{\ll N}u\|_{L_x^{\frac{1}{1-\varepsilon}}L_T^{\frac{2}{1-2\varepsilon}}}^2\\
&\lesssim & \sum_{N\gg 1}\|P_{\ll N}u\|_{L_x^2L_T^{\infty}}^2\|D_x^{\varepsilon}(P_{\gg N}u)^2\|_{L_{xT}^{\frac{2}{1-2\varepsilon}}}^2\\
&\lesssim &T^{1-2\varepsilon}\|u\|_X^4\|P_{\gtrsim 1}u\|_X^2.
\end{eqnarray*}
We repeat the argument of (\ref{eq:D1}), and the proof for (\ref{eq:I_5}) is established.

This concludes the estimate for the contribution of $A_{1,N}$.

\subsection{The contribution of $A_{2,N}$.}
Here the proof is a simple variant of the argument giving (\ref{eq:D1}).
Indeed, the term $P_{\ll N}(i{\cal H}-1)u_xP_{\ll N}uP_+P_Nu$ is quite close to that used for $A_{1,N}$.
Moreover, the proof of Proposition \ref{prop:L^2-spacetime} continues to hold for $i{\cal H}P_{\ll N}u_xP_{\ll N}u$. (Note that for each dyadic number $N\gg 1$, the projection operator ${\cal H}P_N$ is bounded on $L_x^p$ and also on $L_x^pL_T^q$, for $1\le p,q\le \infty$.)

\subsection{The contribution of $A_{3,N}$.}
The proof is a reprise of the argument given in the estimate for $A_{1,N}$.
Note that the integral term $\int_{-\infty}^x{\cal H}P_{\ll N}u_xP_{\ll N}u_x\,dy$ has Fourier support in $|\xi|\ll N$.
In virtue of (\ref{eq:rstrichartz-2}), (\ref{eq:rsmoothing-2}), (\ref{eq:rmaximal-1}), (\ref{eq:rmaximal-2}), we have
\begin{eqnarray}
& & \left(\sum_{N\gg 1}\left\|\int_0^te^{i(t-t')\partial_x^2}A_{3,N}(t')\,dt'\right\|_Y^2\right)^{1/2}
\nonumber\\
& \lesssim & \left(\sum_{N\gg 1}\left\|P_+P_NuD_N(u)\right\|_{L_x^1L_T^2}^2\right)^{1/2}
\label{eq:I-BO4-1}\\
& & +\left(\sum_{N\gg 1}\left(\sum_M\left\|P_M(e^{-\frac{i}{2}\int^x(P_{\ll N}u)^2}P_+P_NuD_N(u))\right\|_{L_x^1L_T^2}\right)^2\right)^{1/2},
\label{eq:I-BO4-2}
\end{eqnarray}
where $D_N(u)=\int_{-\infty}^x{\cal H}P_{\ll N}u_xP_{\ll N}u_x\,dy$.
For the first term (\ref{eq:I-BO4-1}), we split $P_{\ll N}=P_{\lesssim 1}+P_{1\lesssim\cdot\ll N}$ to repeat the argument following the proof of Lemma \ref{lem:L^2}.
In fact, by Lemma \ref{lem:Leibnitz3} and H\"older inequality together with this decomposition, the proof for (\ref{eq:I-BO4-1}) can be reduced to the inequality
\begin{eqnarray*}
\|D_N(u)\|_{L_{xT}^2} & \lesssim & \|D_x^{1/2}P_{\ll N}u\|_{L_{xT}^4}^2\\
&\lesssim & T^{1/2}\|P_{\lesssim 1 }u\|_{L_T^{\infty}L_x^4}^2+\|\langle D_x\rangle^{1/2}u\|_{L_{xT}^4}\|\langle D_x\rangle^{1/2}P_{\gtrsim 1}u\|_{L_{xT}^4},
\end{eqnarray*}
which is bounded by
\begin{eqnarray*}
\lesssim T^{1/2}\|u\|_{X}^2+\|u\|_{X}\|P_{\gtrsim 1}u\|_{X}.
\end{eqnarray*}
On the other hand, the proof for the second term (\ref{eq:I-BO4-2}) follows from combining the above argument with the proof for (\ref{eq:I_1}), which completes the estimate for $A_{3,N}$.

\subsection{The contribution of $A_{4,N}$.}
We may estimate the left-hand side by (\ref{eq:rstrichartz-1}), (\ref{eq:srsmoothing}), (\ref{eq:srmaximal-1}), (\ref{eq:srmaximal-2}).
Therefore it is sufficient to show that
\begin{eqnarray}
& & \left(\sum_{N\gg 1}\left\|P_+P_NuE_N(u)\right\|_{L_{T}^2H_x^{1/2}}^2\right)^{1/2}\label{eq:I-BO5-1}\\
& & +\left(\sum_{N\gg 1}\left(\sum_M\left\|P_M(e^{-\frac{i}{2}\int^x(P_{\ll N}u)^2}P_+P_NuE_N(u))\right\|_{L_T^2H_x^{1/2}}\right)^2\right)^{1/2}
\label{eq:I-BO5-2}\\
&\lesssim & (\|u\|_X^4+\|u\|_X^6)\|P_{\gtrsim 1}u\|_X,
\nonumber
\end{eqnarray}
where $E_N(u)=\int_{-\infty}^xP_{\ll N}(u^2u_x)P_{\ll N}u\,dy$.
We deal with the first term (\ref{eq:I-BO5-1}).
By integrating by parts, observe that
\begin{eqnarray*}
E_N(u) & = & \int_{-\infty}^x\sum_{N_1\sim N_2\ll N}P_{N_1}\partial_xu^3P_{N_2}u\,dy+\int_{-\infty}^x\sum_{N_2\ll N_1\ll N}P_{N_1}\partial_xu^3P_{N_2}u\,dy\\
& & -\int_{-\infty}^x\sum_{N_1\ll N_2\ll N}P_{N_1}u^3P_{N_2}\partial_xu\,dy+\sum_{N_1\ll N_2\ll N}P_{N_1}u^3P_{N_2}u\\.
\end{eqnarray*}
It is easy to bound the contribution of the first term to (\ref{eq:I-BO5-1}) by
\begin{eqnarray*}
& & \lesssim \left(\sum_{N\gg 1}\|P_+P_Nu\|_{L_T^{\infty}H_x^{\frac12}}^2\left\|\sum_{N_1\sim N_2\ll N}\|D_x^{\frac12}P_{N_1}u^3\|_{L_x^2}\|D_x^{\frac12}P_{N_2}u\|_{L_x^2}\right\|_{L_T^2}^2\right)^{1/2}\\
& & \lesssim \|P_{\gg 1}u\|_{X}\|u^3\|_{L_T^{2}H_x^{1/2}}\|u\|_{L_T^{\infty}H_x^{1/2}},
\end{eqnarray*}
which is acceptable since $\|u^3\|_{L_T^2W_x^{1/2,2}}\lesssim\|u\|_{L_{xT}^6}^2\|u\|_{L_T^6W_x^{1/2,6}}\lesssim \|u\|_{X}^3$ (by interpolation).
For the second term , from the Fourier transform, we have
\begin{eqnarray*}
\widehat{\left(\int_{-\infty}^x\sum_{N_2\ll N_1\ll N}P_{N_1}\partial_xu^3P_{N_2}u\,dy\right)}(\xi)=c\sum_{N_2\ll N_1\ll N}\int_{-\infty}^{\infty}\frac{\xi_1}{\xi}\widehat{P_{N_1}u^3}(\xi_1)\widehat{P_{N_2}u}(\xi-\xi_1)\,d\xi_1,
\end{eqnarray*}
so we use the multilinear Fourier multiplier theorem to bound this contribution to (\ref{eq:I-BO5-2}) by
\begin{eqnarray*}
\lesssim\left(\sum_{N\gg 1}\|P_Nu\|_{L_T^4W_x^{1/2,4}}^2\left(\sum_{N_2\ll N_1\ll N}\|P_{N_1}u^3P_{N_2}u\|_{L_{xT}^4}\right)^2\right)^{1/2},
\end{eqnarray*}
which is easy acceptable.
For the third term, the proof is the same as that for the second term.

For the fourth term, we deduce from H\"older inequality that it is
\begin{eqnarray*}
\lesssim \left(\sum_{N\gg 1}\|P_Nu\|_{L_T^{\infty}W_x^{1/2,4}}^2\left(\sum_{N_1\ll N_2\ll N}\|P_{N_1}u^3P_{N_2}u\|_{L_{xT}^4}\right)^2\right)^{1/2},
\end{eqnarray*}
which is acceptable as before.

The proof for (\ref{eq:I-BO5-2}) can be reproduced by combining the above argument with that for (\ref{eq:I_1}).

\subsection{The contribution of $A_{5,N}$.}
The estimate for this contribution is similar but simpler than that for $A_{4,N}$.
This completes the proof for $A_{5,N}$, and hence Proposition \ref{prop:I-BO7}.
\qed
\begin{remark}\label{rem:s>1/2}
We now comment on the case $s>1/2$.
The proof of the above propositions already contains the nonlinear estimates for $s>1/2$.
In particular when $s>1/2$, we require Lemma \ref{lem:retard} with $\theta<1$, in order to use the Leibniz rule on $L_x^pL_T^q$ for $1<p,q<\infty$ (c.f. Lemma \ref{lem:Leibnitz2}).
\end{remark}
\section{Proof of Theorem \ref{thm:k=2}}\label{sec:apriori}
This section is devoted to the proof of Theorem \ref{thm:k=2}.
We shall be concerned with the ``endpoint'' case  $s=1/2$.
(c.f. \cite{MR1} or Remark \ref{rem:s>1/2} for the result for $s>1/2$.)
To begin with, we re-normalize the data a bit via scaling.
By the scaling argument (\ref{eq:scaling}), we have
\begin{eqnarray*}
& & \|u_{0,\lambda}\|_{L^2}=\|u_0\|_{L^2},\\
& & \|u_{0,\lambda}\|_{\dot{H}^{\frac12}}=\frac{1}{\lambda^{1/2}}\|u_0\|_{\dot{H}^{\frac12}}.
\end{eqnarray*}
Thus we may rescale
\begin{eqnarray*}
& & \|P_{\lesssim 1}u_{0,\lambda}\|_{L^2}\le\|u_0\|_{L^2}=C_{low},\\
& & \|(I-P_{\lesssim 1})u_{0,\lambda}\|_{H^{\frac12}}\le\frac{1}{\lambda^{1/2}}\|u_0\|_{H^{\frac12}}<C_{high}\ll 1.
\end{eqnarray*}
Here we choose $\lambda=\lambda(\|u_0\|_{H^{1/2}})\gg 1$, and take the time interval $T$ depending on $\lambda$ later.
We now drop the writing of the scaling parameter $\lambda$ and assume
\begin{eqnarray*}
& & \|P_{\lesssim 1}u_0\|_{L^2}\le C_{low},\\
& & \|(I-P_{\lesssim 1})u_0\|_{H^{\frac12}}\le C_{high}\ll 1.
\end{eqnarray*}
We now apply this to the norms $X$ and $H^{1/2}$, and define new version of the norms of $X$ and $H^{1/2}$, given by with decomposition $I=P_{\lesssim 1}+(I-P_{\lesssim 1})$,
\begin{eqnarray*}
\|u\|_{\widetilde{X}}=\frac{1}{C_{low}}\|P_{\lesssim 1}u\|_X+\frac{1}{C_{high}}\|(I-P_{\lesssim 1})u\|_X, 
\end{eqnarray*}
and
\begin{eqnarray*}
\|\phi\|_{\widetilde{H^{1/2}}}=\frac{1}{C_{low}}\|P_{\lesssim 1}\phi\|_{L^2}+\frac{1}{C_{high}}\|(I-P_{\lesssim 1})\phi\|_{H^{\frac12}}.
\end{eqnarray*}
We remark that $\|u_0\|_{\widetilde{H}^{1/2}}\le 2$.
\subsection{A priori estimate for solutions of (\ref{eq:gBO})}
The purpose of this section is to prove the main a priori estimate for a solution of (\ref{eq:gBO}).
In fact, as a consequence of this estimate, we have the proof of existence, uniqueness and the continuous dependence upon data for the initial value problem (\ref{eq:gBO}).
\begin{proposition}\label{prop:apriori}
Let $u$ be a smooth solution to (\ref{eq:gBO}) and $0<T<1$.
Then we have
\begin{eqnarray}\label{eq:apriori}
\|u\|_{\widetilde{X}}\le C(C_{low})+C(C_{low}+\|u\|_{\widetilde{X}})(T^{\alpha}+C_{high})\|u\|_{\widetilde{X}},
\end{eqnarray}
for some positive $\alpha$.
Here $C(a)\lesssim \langle a\rangle^{100}$.
\end{proposition}
This proposition immediately leads to an a priori estimate for (\ref{eq:gBO}).
\begin{corollary}\label{cor:apriori}
Let $u$ be a smooth solution to (\ref{eq:gBO}).
For $T$ small, $C_{high}$ small, we have
\begin{eqnarray*}
\|u\|_{\widetilde{X}}\le C(C_{high}+C_{low}).
\end{eqnarray*}
\end{corollary}
Before proceeding to the proof of Proposition \ref{prop:apriori}, we establish the following lemmas.
\begin{lemma}\label{lem:low-u}
Let $u$ be a solution to (\ref{eq:gBO}).
Then
\begin{eqnarray*}
\|P_{\lesssim 1}u\|_{X}\lesssim C_{low}+T^{1/2}\|u\|_{X}^3.
\end{eqnarray*}
\end{lemma}
{\it Proof of Lemma \ref{lem:low-u}.}
Applying $P_+$ to (\ref{eq:gBO}), we obtain the equation
\begin{eqnarray*}
(\partial_t-i\partial_x^2)P_+u=P_+(u^2u_x).
\end{eqnarray*}
Using the integral equation
\begin{eqnarray*}
P_+u(t)=e^{it\partial_x^2}P_+u_0-\int_0^te^{i(t-t')\partial_x^2}P_+(u^2u_x)(t')\,dt',
\end{eqnarray*}
and by (\ref{eq:strichartz}), (\ref{eq:smoothing}), (\ref{eq:maximal}), (\ref{eq:smoothing-1}), (\ref{eq:rstrichartz-1}), (\ref{eq:srsmoothing}), (\ref{eq:srmaximal-1}), (\ref{eq:srmaximal-2}), we have
\begin{eqnarray*}
\|P_{\lesssim 1}P_+u\|_{X} & \lesssim & \|P_{\lesssim 1}P_+u_0\|_{H^{1/2}}+\|P_{\lesssim 1}P_+(u^3)_x\|_{L_T^1H_x^{1/2}}\\
& \lesssim & C_{low}+T^{1/2}\|u\|_{L_{xT}^6}^2\|u\|_{L_T^6W_x^{1/2,6}}.
\end{eqnarray*}
Since $u$ is real-valued, this proves Lemma \ref{lem:low-u}.
\qed
\begin{lemma}\label{lem:u-v}
Let $u$ and $v_N$ be given in (\ref{eq:gauge}).
Then
\begin{eqnarray*}
\|P_{\gg 1}u\|_{X} \lesssim (1+\|u\|_{L_T^{\infty}H^{1/2}}^4)\left(\sum_{N\gg 1}\|v_N\|_{Y}^2\right)^{\frac12},
\end{eqnarray*}
where the space $Y$ is defined in section \ref{sec:nonlinear}.
\end{lemma}
{\it Proof of Lemma \ref{lem:u-v}.}
We will consider separately each of contribution of $L_T^{\infty}H^{1/2}$, $L_x^{\infty}L_T^2$, $L_x^2L_T^{\infty}$  and $L_x^4L_T^{\infty}$-norms.

To bound the contribution of the $L_T^{\infty}H^{1/2}$-norm, since $u$ is real,  we use Leibniz' rule (c.f. Lemma \ref{lem:Leibnitz2}) to estimate 
\begin{eqnarray*}
\|D_x^{\frac12}P_N(e^{\frac{i}{2}\int^x(P_{\ll N}u)^2}v_N)\|_{L^2}\lesssim (1+\|P_{\ll N}u\|_{H^{\frac12}}^2)\|v_N\|_{H^{\frac12}},
\end{eqnarray*}
which gives the estimate, summing on $l_N^2$.

For the contribution of the $L_x^{\infty}L_T^2$-norm, observe first that
\begin{eqnarray*}
\partial_xP_+P_Nu=e^{\frac{i}{2}\int^x(P_{\ll N}u)^2}\left(\partial_xv_N+\frac{i}{2}(P_{\ll N}u)^2v_N\right),
\end{eqnarray*}
then by Littlewood-Paley decomposition, we bound the $L_x^{\infty}L_T^2$- norm of $\partial_xP_Nu$ by
\begin{eqnarray}\label{eq:smo}
\|\partial_xP_Nu\|_{L_x^{\infty}L_T^2}\lesssim \|\partial_xv_N\|_{L_x^{\infty}L_T^2}+\left\|\widetilde{P}_{N}\left(\sum_{N_1}P_{N_1}(e^{\frac{i}{2}\int^x(P_{\ll N}u)^2})(P_{\ll N}u)^2\sum_{N_2}P_{N_2}v_N\right)\right\|_{L_x^{\infty}L_T^2}.
\end{eqnarray}
To estimate the second term, we first split $\sum_{N_2}=\sum_{N_2\sim N}+\sum_{N_2\not\sim N}$.
For $N_2\sim N$, we bound this contribution to the second term of (\ref{eq:smo}) by
\begin{eqnarray*}
c\|(P_{\ll N}u)^2\sum_{N_2\sim N}P_{N_2}v_N\|_{L_x^{\infty}L_T^2}& \lesssim & \|P_{\ll N}u\|_{L_{xT}^{\infty}}^2\sum_{N_2\sim N}\|P_{N_2}v_N\|_{L_x^{\infty}L_T^2}\\
& \lesssim & \|u\|_{L_T^{\infty}H^{1/2}}^2\sum_{N_2\sim N}\|P_{N_2}\partial_xv_N\|_{L_x^{\infty}L_T^2},
\end{eqnarray*}
which is acceptable.

In $N_2\not\sim N$, we split again $\sum_{N_2\not\sim N}=\sum_{N_2\ll N}+\sum_{N_2\gg N}$.
For $N_2\ll N$, observe that
\begin{eqnarray*}
\widetilde{P}_N(\sum_{N_1}P_{N_1}e^{\frac{i}{2}\int^x(P_{\ll N}u)^2}(P_{\ll N}u)^2P_{\ll N}v_N)=\widetilde{P}_N(P_{\sim N}e^{\frac{i}{2}\int^x(P_{\ll N}u)^2}(P_{\ll N}u)^2P_{\ll N}v_N),
\end{eqnarray*}
while for $N_2\gg N$, we see that the left hand side
\begin{eqnarray*}
=\widetilde{P}_N(\sum_{N_1\sim N_2\gg N}P_{N_1}e^{\frac{i}{2}\int^x(P_{\ll N}u)^2}(P_{\ll N}u)^2P_{N_2}v_N).
\end{eqnarray*}
Then we have the bound of this contribution to the second term of (\ref{eq:smo}) by
\begin{eqnarray*}
& & \lesssim \|P_{\sim N}e^{\frac{i}{2}\int^x(P_{\ll N}u)^2}\|_{L_{xT}^{\infty}}\|P_{\ll N}u\|_{L_{xT}^{\infty}}^2\|P_{\ll N}v_N\|_{L_{xT}^{\infty}}\\
& & +\sum_{N_1\sim N_2\gg N}\|P_{N_1}e^{\frac{i}{2}\int^x(P_{\ll N}u)^2}\|_{L_{xT}^{\infty}}\|P_{\ll N}u\|_{L_{xT}^{\infty}}^2\|P_{N_2}v_N\|_{L_{xT}^{\infty}}\\
& & \lesssim  \|u\|_{L_T^{\infty}H^{\frac12}}^4\|v_N\|_{L_T^{\infty}H^{\frac12}}.
\end{eqnarray*}
Therefore summing also on $l_N^2$, we complete the proof for the contribution of the $L_x^{\infty}L_T^2$-norm.

The estimate for the contribution of the $L_x^2L_T^{\infty}$-norm is easy, since $|P_+P_Nu|=|v_N|$.

The proof for the contribution of the $L_x^4L_T^{\infty}$-norm is in the same style as that for the $L_x^{\infty}L_T^2$-norm, because $\|\langle D_x\rangle^{1/4}P_Nu\|_{L_x^4L_T^{\infty}}\sim N^{1/4}\|P_Nu\|_{L_x^4L_T^{\infty}}$.
We reprise the argument following the proof for the contribution of the $L_x^{\infty}L_T^2$-norm, to obtain the bound
\begin{eqnarray*}
& & \|\langle D_x\rangle^{1/4}P_Nu\|_{L_x^4L_T^{\infty}}\\
&\lesssim & N^{\frac14}\sum_{N_2\sim N}\|P_{N_2}v_N\|_{L_x^4L_T^{\infty}}+N^{1/4}\|P_{\sim N}e^{\frac{i}{2}\int^x(P_{\ll N}u)^2}(P_{\ll N}u)^2P_{\ll N}v_N\|_{L_x^4L_T^{\infty}}\\
& & +\sum_{N_1\sim N_2\gg N}N^{1/4}\|P_{N_1}e^{\frac{i}{2}\int^x(P_{\ll N}u)^2}(P_{\ll N}u)^2P_{N_2}v_N\|_{L_x^4L_T^{\infty}}\\
& \lesssim & (1+\|u\|_{L_T^{\infty}H^{1/2}}^4)\|\langle D_x\rangle^{\frac14}v_N\|_{L_x^4L_T^{\infty}}.
\end{eqnarray*}
We apply $l_N^2$-sum and thus prove the estimate for the contribution of the $L_x^4L_T^{\infty}$-norm.

This concludes the proof of Lemma \ref{lem:u-v}.
\qed

For the proof of Proposition \ref{prop:apriori}, we will use the following estimate concerning (\ref{eq:I-BO1}). 
\begin{lemma}\label{lem:linear-v}
For $\phi\in H^{1/2}$,
\begin{eqnarray*}
\left(\sum_{N\gg 1}\|e^{it\partial_x^2}(e^{-\frac{i}{2}\int^x(P_{\ll N}\phi)^2}P_+P_N\phi)\|_{Y}^2\right)^{1/2}\lesssim (1+\|\phi\|_{H^{\frac12}}^2)\|P_{\gg 1}\phi\|_{H^{\frac12}}.
\end{eqnarray*}
\end{lemma}
{\it Proof of Lemma \ref{lem:linear-v}.}
Applying (\ref{eq:strichartz}), (\ref{eq:smoothing}), (\ref{eq:maximal}) and (\ref{eq:smoothing-1}) shows that it is sufficient to prove
\begin{eqnarray}
& & \left(\sum_{N\gg 1}\|e^{-\frac{i}{2}\int^x(P_{\ll N}\phi)^2}P_+P_N\phi\|_{H^{\frac12}}^2\right)^{\frac12}+\left(\sum_{N\gg 1}\left(\sum_M\|P_M(e^{-\frac{i}{2}\int^x(P_{\ll N}\phi)^2}P_+P_N\phi)\|_{H^{\frac12}}\right)^2\right)^{\frac12}
\nonumber\\
& & \lesssim  (1+\|\phi\|_{H^{\frac12}}^2)\|P_{\gg 1}\phi\|_{H^{\frac12}}.
\label{eq:li}
\end{eqnarray}

By Leibniz' rule (c.f. Lemma \ref{lem:Leibnitz2}), the first term of the left hand side of (\ref{eq:li}) is bounded by
\begin{eqnarray*}
\lesssim (1+\|\phi\|_{H^{\frac12}}^2)\|P_{\gg 1}\phi\|_{H^{\frac12}},
\end{eqnarray*}
which leads to a desired estimate.

Next, we deal with the second term.
Like the argument in the proof in section \ref{sec:L^2-spacetime}, we split $\sum_M=\sum_{M\lesssim N}+\sum_{M\gg N}$.
Hence we bound this contribution to the left-hand side of (\ref{eq:li}) by
\begin{eqnarray*}
& & \lesssim \left(\sum_{N\gg 1}\|P_N\phi\|_{H^{\frac12}}^2\right)^{1/2}+\left(\sum_{N\gg 1}\left(\sum_{M\gg N}\|\widetilde{P}_Me^{-\frac{i}{2}\int^x(P_{\ll N}\phi)^2}\|_{H^{\frac12}}\right)^2\|P_N\phi\|_{L^{\infty}}^2\right)^{\frac12}\\
& \lesssim & (1+\|\phi\|_{H^{\frac12}}^2)\|P_{\gg 1}\phi\|_{H^{\frac12}},
\end{eqnarray*}
which is also acceptable.
\qed
\newline
{\it Proof of Proposition \ref{prop:apriori}.}

Turning to the proof of Proposition \ref{prop:apriori}, with the above lemmas, we show the a priori estimate for solutions of (\ref{eq:gBO}).

In light of Lemma \ref{lem:u-v}, it is reasonable to pass from the a priori estimate for $u$ to that of $v_N$.
We deduce from Proposition \ref{prop:I-BO7} 
together with Lemmas \ref{lem:low-u}, \ref{lem:u-v} and \ref{lem:linear-v} that
\begin{eqnarray*}
\|u\|_{\widetilde{X}} & \le & C(C_{low})+C(C_{low})T^{\frac12}\|u\|_{\widetilde{X}}^3\\
& & +C(C_{low})(1+\|u\|_{L_T^{\infty}H_x^{1/2}}^4)\left(1+T^{\frac14}(1+\|u\|_{\widetilde{X}}^6)\|u\|_{\widetilde{X}}+C_{high}^2\|u\|_{\widetilde{X}}+C_{high}(1+\|u\|_{\widetilde{X}}^2)\|u\|_{\widetilde{X}}^2\right).
\end{eqnarray*}
Observe that by renormalization of $\|\cdot\|_{H^{1/2}}$-norm we see that
\begin{eqnarray*}
\|u(t)\|_{H^{1/2}}\lesssim \|P_{\lesssim 1}u(t)\|_{L^2}
+C_{high}\|(I-P_{\lesssim 1})u(t)\|_{\widetilde{H}^{1/2}}.
\end{eqnarray*}
The high frequency part $C_{high}\|(I-P_{\lesssim 1})u\|_{\widetilde{H}_x^{1/2}}$ can be absorbed into the $\widetilde{X}$-norm.
Then substituting Lemma \ref{lem:low-u} again in estimating the low frequency part of the norm $\|P_{\lesssim 1}u\|_{L_T^{\infty}H_x^{1/2}}$, we obtain (\ref{eq:apriori}) and complete the proof of Proposition \ref{prop:apriori}.
\qed
\subsection{Proof of Theorem \ref{thm:k=2}.}
We come now to the proof of Theorem \ref{thm:k=2}, and describe the key points when we follow the compactness argument with the a priori estimate.
We refer to the papers \cite{MR1,MR2,KT,P,Ta1} for the details. 

Let $\{u_{0n}\}$ be a sequence in $H^{\infty}$ such that $u_{0n}\to u_0$ in $H^{1/2}$ as $n\to\infty$ and, $\|u_{0n}\|_{H^{1/2}}\le 2\|u_0\|_{H^{1/2}}$.
We see that if $u_n$ is a $H^{\infty}$-solution to (\ref{eq:gBO}) with data $u_n(0)=u_{0n}$, then we have the a priori estimate (\ref{eq:apriori}): with Corollary \ref{cor:apriori}, for small $T>0$ (we take $C_{high}$ small), 
\begin{eqnarray}\label{eq:n1-apriori}
\|u_n\|_{\widetilde{X}}\lesssim C(\|u_0\|_{H^{1/2}}).
\end{eqnarray}
Similarly, noting that $|e^{i\int^xf_1}-e^{i\int^xf_2}|\lesssim \|f_1-f_2\|_{L^1}$ for real functions $f_1,f_2$, we obtain
\begin{eqnarray}\label{eq:n2-apriori}
\|u_n-u_{n'}\|_{\widetilde{X}}\lesssim C(\|u_0\|_{H^{1/2}})\|u_{0n}-u_{0n'}\|_{\widetilde{H}^{1/2}}
\end{eqnarray}
(by using estimates similar to (\ref{eq:apriori}) for differences of solutions).
These bounds (\ref{eq:n1-apriori}) and (\ref{eq:n2-apriori}) will allow us to obtain the existence of the solution $u\in\widetilde{X}$ to (\ref{eq:gBO}).
In particular, using Fatou's lemma, we can show
\begin{eqnarray*}
\|u\|_{\widetilde{X}}\lesssim C(\|u_0\|_{H^{1/2}}).
\end{eqnarray*}

Now we prove the uniqueness of solution.
Let $u$ and $\widetilde{u}$ be two solutions of (\ref{eq:gBO}) with data $u_0$ and $\widetilde{u}_0$, respectively.
By (\ref{eq:n2-apriori}) (choose $T>0$ and $C_{high}$ smaller, if necessary), we have
\begin{eqnarray*}
\|u-\widetilde{u}\|_{\widetilde{X}}\lesssim C_0(\|u_0\|_{H^{1/2}}+\|\widetilde{u}_0\|_{H^{1/2}})\|u_0-\widetilde{u}_0\|_{\widetilde{H}^{1/2}}.
\end{eqnarray*}
Thus the solution is unique in $\widetilde{X}$, also in $X$.

The continuous dependence of solution on data is actually proven in the same way as in the proof of the existence of solution.

This concludes the proof of Theorem \ref{thm:k=2}.
\qed

Carlos E. Kenig, Department of Mathematics,
                 University of Chicago
                 Chicago, Illinois 60637, USA

Hideo Takaoka, Department of Mathematics,
               Kobe University,
               Kobe 657-8501, Japan
\end{document}